\pdfoutput=1
\RequirePackage{ifpdf}
\ifpdf % We~are running pdfTeX in pdf mode
\documentclass[pdftex]{sigma}
\else
\documentclass{sigma}
\fi

\def \bui#1#2{\mathrel{\mathop{\kern 0pt#1}\limits^{#2}}}
\def \buil#1#2{\mathrel{\mathop{\kern 0pt#1}\limits_{#2}}}

\let\<\langle
\let\>\rangle
\newcommand{\R}{{\mathbb R}}

\numberwithin{equation}{section}

\newtheorem{Theorem}{Theorem}[section]
\newtheorem*{Theorem*}{Theorem}
\newtheorem{Corollary}[Theorem]{Corollary}
\newtheorem{Lemma}[Theorem]{Lemma}
\newtheorem{Proposition}[Theorem]{Proposition}

\theoremstyle{definition}

\newtheorem{Remark}[Theorem]{Remark}
\newtheorem{Remarks}[Theorem]{Remarks}

\begin{document}

\allowdisplaybreaks

\newcommand{\arXivNumber}{2412.05612}

\renewcommand{\PaperNumber}{010}

\FirstPageHeading

\ShortArticleName{The Buckling and Clamped Plate Problems on Differential Forms}

\ArticleName{The Buckling and Clamped Plate Problems\\ on Differential Forms}

\Author{Fida EL CHAMI~$^{\rm a}$, Nicolas GINOUX~$^{\rm b}$, Georges HABIB~$^{\rm ab}$, Ola MAKHOUL~$^{\rm a}$\newline and Simon RAULOT~$^{\rm c}$}

\AuthorNameForHeading{F.~El Chami, N.~Ginoux, G.~Habib, O.~Makhoul and S.~Raulot}

\Address{$^{\rm a)}$~Department of Mathematics, Faculty of Sciences II, Lebanese University,\\
\hphantom{$^{\rm a)}$}~P.O. Box 90656 Fanar-Matn, Lebanon}
\EmailD{\mail{fchami@ul.edu.lb}, \mail{ghabib@ul.edu.lb}, \mail{ola.makhoul@ul.edu.lb}}

\Address{$^{\rm b)}$~Universit\'e de Lorraine, CNRS, IECL, 57000 Metz, France}
\EmailD{\mail{nicolas.ginoux@univ-lorraine.fr}}
\URLaddressD{\url{https://nicolas-ginoux.perso.math.cnrs.fr/}}

\Address{$^{\rm c)}$~Universit\'e de Rouen Normandie, CNRS, Normandie Univ, LMRS UMR 6085,\\
\hphantom{$^{\rm c)}$}~76000 Rouen, France}
\EmailD{\mail{simon.raulot@univ-rouen.fr}}
\URLaddressD{\url{https://lmrs.univ-rouen.fr/persopage/simon-raulot}}

\ArticleDates{Received July 17, 2025, in final form January 22, 2026; Published online February 03, 2026}

\Abstract{We extend the buckling and clamped-plate problems to the context of differential forms on compact Riemannian manifolds with smooth boundary. We characterize their smallest eigenvalues and prove that, in the case of bounded Euclidean domains, their spectra without multiplicities on forms coincide with the spectra of the corresponding problems on functions.
We obtain various estimates involving the first eigenvalues of the mentioned problems and the ones of the Hodge Laplacian with respect to Dirichlet and absolute boundary conditions on forms.
These estimates generalize previous ones in the case of functions.}


\Classification{58J32; 58C40}

\section{Introduction}
Let $(M,g)$ be an $n$-dimensional compact Riemannian manifold with smooth boundary $\partial M$ and let $\nu$ be the inward unit vector field normal to $\partial M$.
For a smooth function $f$ on $M$, we consider the following two problems:
\begin{alignat}{3}
&\Delta^2 f = \Gamma f\qquad && \textrm{on } M,&\nonumber\\
&f = 0 \quad && \textrm{on } \partial M,&\label{clamped_plate_functions}\\
&\frac{\partial f}{\partial\nu} = 0 \quad&& \textrm{on } \partial M&\nonumber
\end{alignat}
and
\begin{alignat}{3}
&\Delta^2 f=\Lambda \Delta f \qquad && \textrm{on } M,&\nonumber\\
&f = 0 \quad && \textrm{on } \partial M,&\label{buckling_functions}\\
&\frac{\partial f}{\partial\nu} = 0 \quad && \textrm{on } \partial M&\nonumber
\end{alignat}
called the clamped plate and the buckling problem respectively.
Note that \smash{$\Delta f=- \hbox{tr}\bigl(\nabla^2f\bigr)$} is the Laplace operator of $f$ and $\Delta^2 $ its square, which is sometimes called the bi-Laplace operator.
It is well known that these two problems have discrete spectra consisting of eigenvalues of finite multiplicities
\[ 0< \Gamma_1\le \Gamma_2 \le \dots \le \Gamma_k \le \dots \longrightarrow \infty\]
and
\[ 0< \Lambda_1 \le \Lambda_2 \le \dots \le \Lambda_k \le \dots \longrightarrow \infty,\]
where each eigenvalue is repeated according to its multiplicity.
Physically, problem $(\ref{clamped_plate_functions})$ describes the vibrations of a clamped plate, whereas problem $(\ref{buckling_functions})$ describes the critical buckling load of a clamped plate subjected to a uniform compressive force around its boundary.

These two problems were studied by numerous authors.
In 1955, Payne~\cite{Payne1955} proved that, if~$M$ is a planar bounded domain, then
\begin{equation*}
\Lambda_1 \ge \lambda_2,
\end{equation*}
where $\lambda_2$ is the second eigenvalue of the Dirichlet problem on $M$ (see also~\cite{Friedlander2004} for a corrected proof).
In 1996, Ashbaugh and Laugesen showed, in their work in~\cite{AL}, that whenever $M$ is a~bounded and connected open subset of the Euclidean space $\R^n$,
\begin{equation}
\label{AL} \Lambda_1^2 \ge \Gamma_1 \ge \Lambda_1 \lambda_1 > \lambda_1^2,
\end{equation}
where $\lambda_1$ is the first eigenvalue of the Dirichlet problem on $M$.
 In~\cite{ChenChengWangXia2012}, Chen, Cheng, Wang and Xia proved that, if the Ricci curvature is bounded below by $n-1$, then
\begin{equation*}%\label{CCWX}
\Gamma_1 > n\lambda_1\qquad\hbox{and}\qquad \Lambda_1>n.
\end{equation*}
Ilias and Shouman gave in~\cite{IliasShouman2019} an estimate relating the first eigenvalue $\Lambda_1$ of the buckling problem to $\mu_1$, the first nonzero eigenvalue of the Neumann problem
 \begin{equation} \label{IS}
 \mu_1 < \Lambda_1.
 \end{equation}

In this paper, we first generalize problems~(\ref{clamped_plate_functions}) and~(\ref{buckling_functions}) to the context of differential forms on the manifold $M$.
We prove that each problem has a discrete spectrum consisting of a nondecreasing sequence of real eigenvalues of finite multiplicities and the corresponding eigenforms form a Hilbert basis of $L^2$-integrable $p$-forms on $M$.
We also characterize the first eigenvalue of each problem (see Theorems~\ref{BucKPb} and~\ref{CP_Pb}).

In Section~\ref{s:ebs}, we prove that if $(M,g)$ is a bounded domain of the Euclidean space $\R^n$, the spectra without multiplicities of both problems on $p$-forms, for $p=1,\dots,n$, and on functions on~$M$ coincide (see Proposition~\ref{BCP_Euclidean}).
This allows us, for example, to determine the first eigenvalues of both of the problems on $p$-forms for the Euclidean ball of arbitrary radius.
In the same section, namely in Theorem~\ref{BuckCP}, we establish a relationship between the first eigenvalues of the buckling and clamped-plate problems on an arbitrary compact Riemannian manifold $(M,g)$ with smooth boundary.
In the same context, we consider the following problem, called the Dirichlet problem on forms,
\begin{equation}
\begin{aligned}
&\Delta \omega =\lambda \omega\ \quad &&\textrm{in } M, \\
&\omega = 0\ \quad &&\textrm{on } \partial M
\end{aligned}\label{Dirichlet-forms}
\end{equation}
and give estimates relating the first eigenvalues of the buckling, clamped-plate and Dirichlet problems. The estimates that we obtain generalize inequalities~(\ref{AL}) to differential forms on~$M$.
We prove as well that there exists a connection %between the first eigenvalue of the buckling problem on forms
with the first eigenvalue of the absolute boundary value problem on forms
\[
\begin{aligned}
&\Delta \omega =\mu \omega \ \quad&&\textrm{on } M, \\
&\nu\lrcorner \omega = 0 \ \quad&& \textrm{on } \partial M, \\
&\nu \lrcorner \, {\rm d}\omega = 0 \ \quad&& \textrm{on } \partial M
\end{aligned}
\]
given in Theorem~\ref{BuckAbs}.
This connection allows to extend~(\ref{IS}) to differential forms.

Further in the same section, we show that, under the condition that the Weitzenb\"ock curvature operator is bounded below by a positive constant $\gamma$, the first eigenvalue $\lambda_{1,p}$ of problem~(\ref{Dirichlet-forms}) is also bounded below by a quantity depending on $\gamma$, see Theorem~\ref{GM_Dirichlet}.
This gives new estimates of the first eigenvalues of the buckling and clamped-plate problems under the same conditions, see Corollary~\ref{BCDL}.

We end Section~\ref{s:ebs} by considering domains $M$ of the unit round sphere $\mathbb{S}^n$ and derive inequalities relating the first eigenvalues of the buckling and clamped plate problems on forms of different degrees (see Theorem~\ref{t:domainSn}).

Finally, we point out that our boundary conditions for the buckling and clamped plate problems are the most natural in the context of differential forms.
To what extent there exist further boundary conditions generalizing the ones for the scalar problem is a question which has not been addressed in this article.
Independently, so-called universal inequalities have been established in the scalar case~\cite{ChengYang2006,ChengYang2012,DuMaoWang2016}.
Extending those to our problems constitutes the object of future work.

\section{The buckling and clamped plate problems}

In order to make this work self-contained, we collect here some classical and useful formulae in the study of $p$-forms on manifolds with boundary. In the following, $(M^n,g)$ denotes an $n$-dimensional compact Riemannian manifold with smooth boundary and $\nu$ the inner unit normal to $\partial M$. Then we recall from~\cite[Lemma 18]{RaulotSavo2011} that for any $p$-form $\omega$ on $M$, we have
\begin{equation}\label{eq:iotastarnablanuomega}
\iota^*(\nabla_\nu\omega)=\nu\lrcorner\,{\rm d}\omega+{\rm d}(\nu\lrcorner\,\omega)+S^{[p]}(\iota^*\omega)
\end{equation}
and
\begin{equation}\label{eq:nuintnablanuomega}
\nu\lrcorner\,\nabla_\nu\omega = \delta^{\partial M}(\iota^*\omega)-\iota^*(\delta\omega)-S^{[p-1]}(\nu\lrcorner\,\omega)+(n-1)H\nu\lrcorner\,\omega.
\end{equation}
Here \smash{$\delta^{\partial M}$} denotes the codifferential on \smash{$\partial M$}, \smash{$S^{[p]}$} is the natural extension as an endomorphism of $\Omega^{p}(\partial M)$ of the shape operator $S:=-\nabla\nu$ of the embedding $\iota$ of $\partial M$ in $M$ and \smash{$H=\frac{1}{n-1}\operatorname{tr}(S)$} is its mean curvature, see, e.g.,~\cite[p.~624]{RaulotSavo2011}.

We are especially interested in $p$-forms $\omega\in\Omega^p(M)$ which satisfy the boundary condition $\omega_{|_{\partial M}}=0$.
Then, using~(\ref{eq:iotastarnablanuomega}) and~(\ref{eq:nuintnablanuomega}), it is not difficult to compute that, along $\partial M$, we have $\iota^*\nabla_\nu\omega=\nu\lrcorner\,{\rm d}\omega$ and $\nu\lrcorner\nabla_\nu\omega=-\iota^*\delta\omega$.
So it is straightforward to see that
\begin{equation}\label{Zcharacterization}
\begin{aligned}
&\omega_{|\partial M} = 0,\\
&\nabla_\nu\omega_{|\partial M} = 0,
\end{aligned}
\quad\Longleftrightarrow\quad
\begin{aligned}
& \iota^*\omega = 0,\\
& \nu\lrcorner\omega = 0,\\
& \iota^*\delta\omega = 0,\\
& \nu\lrcorner\,{\rm d}\omega = 0.
\end{aligned}
\end{equation}
In the following, we will denote by $Z$ the vector space of smooth $p$-forms on $M$ which satisfy these boundary conditions that is
\begin{equation*}%\label{DefZ}
Z:=\big\{\omega\in\Omega^p(M)\mid \omega_{|\partial M} = 0\text{ and }\nabla_\nu\omega_{|\partial M} = 0\big\}.
\end{equation*}

In this work, we will also often use integration by parts formulae which are cumbersome to write down in the general framework of $p$-forms on manifolds with boundary. However, as we restrict our attention to elements in $Z$, it may be very useful to observe that in this context they become very simple.
In fact, for any \smash{$\omega,\omega'\in\Omega^p(M)$}, it holds that
\begin{gather}
\int_M\bigl\langle\Delta\omega,\omega'\bigr\rangle {\rm d}\mu = \int_M\bigl(\bigl\langle {\rm d}\omega,{\rm d}\omega'\bigr\rangle +\bigl\langle\delta\omega,\delta\omega'\bigr\rangle\bigr) {\rm d}\mu \nonumber\\
\hphantom{\int_M\bigl\langle\Delta\omega,\omega'\bigr\rangle {\rm d}\mu =}{}
+\int_{\partial M}\bigl(\bigl\langle \nu\lrcorner {\rm d}\omega,\iota^*\omega'\bigr\rangle-\bigl\langle\iota^*\delta\omega,\nu\lrcorner\omega'\bigr\rangle\bigr) {\rm d}\sigma \label{eq:laplacianG1}\\
\hphantom{\int_M\bigl\langle\Delta\omega,\omega'\bigr\rangle {\rm d}\mu}{}
= \int_M\bigl\langle\omega,\Delta\omega'\bigr\rangle {\rm d}\mu + \int_{\partial M}\bigl(\bigl\langle\nu\lrcorner
{\rm d}\omega,\iota^*\omega'\bigr\rangle-\bigl\langle\iota^*\omega,\nu\lrcorner
{\rm d}\omega'\bigr\rangle
\nonumber\\
\hphantom{\int_M\bigl\langle\Delta\omega,\omega'\bigr\rangle {\rm d}\mu =\int_M\bigl\langle\omega,\Delta\omega'\bigr\rangle {\rm d}\mu + \int_{\partial M}}{}
+\bigl\langle\nu\lrcorner\omega,
\iota^*\delta\omega'\bigr\rangle-\bigl\langle\iota^*\delta\omega,
\nu\lrcorner\omega'\bigr\rangle\bigr){\rm d}\sigma,\label{eq:laplacianG}
\end{gather}
and so we immediately deduce from~(\ref{Zcharacterization}) that if $\omega,\omega'\in Z$, then
\begin{equation}\label{eq:laplacian}
\int_M\bigl\langle\Delta\omega,\omega'\bigr\rangle {\rm d}\mu=\int_M\bigl(\bigl\langle {\rm d}\omega,{\rm d}\omega'\bigr\rangle+\bigl\langle\delta\omega,\delta\omega'\bigr\rangle\bigr){\rm d}\mu=\int_M\bigl\langle\omega,\Delta\omega'\bigr\rangle {\rm d}\mu.
\end{equation}
Here ${\rm d}\mu$ (resp.\ ${\rm d}\sigma$) denotes the Riemmanian measure density of $(M^n,g)$ (resp.\ $\partial M$ endowed with the induced metric). Now by replacing $\omega$ by $\Delta \omega$ in~(\ref{eq:laplacianG}), we obtain
\begin{align}
\nonumber\int_M\bigl\langle\Delta^2\omega,\omega'\bigr\rangle {\rm d}\mu ={}& \int_M\bigl\langle\Delta \omega,\Delta\omega'\bigr\rangle {\rm d}\mu+\int_{\partial M}\bigl(\bigl\langle\nu\lrcorner {\rm d}\Delta\omega,\iota^*\omega'\bigr\rangle-\bigl\langle\iota^*\Delta\omega,\nu\lrcorner {\rm d}\omega'\bigr\rangle\bigr){\rm d}\sigma\\
 &{}{+}\,\int_{\partial M}\bigl(\bigl\langle\nu\lrcorner\Delta\omega,\iota^*\delta\omega'\bigr\rangle-\bigl\langle\iota^*\delta\Delta\omega,\nu\lrcorner\omega'\bigr\rangle\bigr){\rm d}\sigma\label{eq:partialintegDelta2omega}
\end{align}
and so if $\omega,\omega'\in Z$, we get
\begin{equation}\label{eq:partialintDeltacons}
\int_M\bigl\langle\Delta^2\omega,\omega'\bigr\rangle {\rm d}\mu=\int_M\bigl\langle\Delta \omega,\Delta\omega'\bigr\rangle {\rm d}\mu.
\end{equation}
In the following, we will also denote by $(\cdot,\cdot)_{L^2(M)}$ the $L^2$-scalar product on $\Omega^p(M)$ and $\|\cdot\|_{L^2(M)}$ its associated norm.
We finally notice that the boundary conditions studied here turn out to be elliptic in the sense of Lopatinski\u{\i}--Shapiro (see~\cite[Definition~1.6.1]{Sc}). This was proved in a~more general setting by the first four authors (see~\cite[Lemma~6.1]{EGHM}) and we restate this result in our context.
\begin{Lemma}\label{elliptic_problem}
Let $(M^n,g)$ be a compact Riemannian manifold with smooth boundary $\partial M$ and let~$\nu$ be the inward unit normal vector field to the boundary. The following boundary value problem:
\[
\begin{aligned}
&\Delta^2\omega=f\ \quad &&\textrm{on }M,\\
&\omega = \omega_1 \ \quad && \textrm{on }\partial M,\\
&\iota^*\delta\omega =\omega_2\ \quad &&\textrm{on }\partial M,\\
&\nu\lrcorner {\rm d}\omega =\omega_3\ \quad && \textrm{on }\partial M
\end{aligned}
\]
for given $f\in\Omega^p(M)$, $\omega_1\in\Omega^p(M)_{|_{\partial M}}$, $\omega_2\in\Omega^{p-1}(\partial M)$, $\omega_3\in\Omega^p(\partial M)$, is elliptic in the sense of Lopatinski\u{\i}--Shapiro.
\end{Lemma}

\subsection{The buckling problem}

The buckling eigenvalue problem on differential forms is
\begin{equation}\label{buckling_forms}
\begin{aligned}
&\Delta^2 \omega=\Lambda \Delta \omega \ \quad &&\textrm{on } M, \\
&\omega=0\ \quad &&\textrm{on }\partial M, \\
&\nabla_\nu\omega_{|\partial M}=0\ \quad &&\textrm{on }\partial M ,
\end{aligned}
\end{equation}
for some real constant $\Lambda$.
Let us begin with the following existence result.
\begin{Theorem}\label{BucKPb}
There exists a Hilbert basis of the space of $L^2$-integrable $p$-forms on $(M^n,g)$ consisting of eigenforms solutions of the problem~\eqref{buckling_forms} associated to an unbounded and positive sequence of eigenvalues $(\Lambda_{i,p})_{i\geq 1}$. Moreover, each eigenspace has a finite multiplicity and the corresponding eigenforms are smooth.
Finally, the first eigenvalue $\Lambda_{1,p}$ is characterized by
\begin{gather}\label{eq:charb}
\Lambda_{1,p}=\inf\Biggl\{\frac{\|\Delta\omega\|_{L^2(M)}^2}{
\| {\rm d}\omega\|_{L^2(M)}^2+\|\delta \omega \|_{L^2(M)}^2},\, \omega\in\Omega^p(M)\setminus\{0\},\, \omega_{|_{\partial
M}}=0,\, \nabla_\nu\omega_{|\partial M}=0\Biggr\}.
\end{gather}
Equality holds if and only if $\omega$ is an eigenform associated to the first eigenvalue.
\end{Theorem}

\begin{proof}
The proof is classical and so we only recall the main steps. For all $p$-forms $\omega$, $\omega'$, the two bilinear forms
\[
\bigl(\omega,\omega'\bigr)_V:=\int_M\bigl\langle\Delta\omega,\Delta\omega'\bigr\rangle
{\rm d}\mu \qquad \textrm{and}\qquad \bigl(\omega,\omega'\bigr)_W:=\int_{M}\bigl(
\bigl\langle {\rm d}\omega,{\rm d}\omega'\bigr\rangle+\bigl\langle \delta\omega,\delta\omega'\bigr\rangle\bigr) {\rm d}\mu
\]
define scalar products on $Z$ whose associated norms will be denoted by $\|\cdot\|_V$ and $\|\cdot\|_W$.
We~will also denote by $V$ and $W$ the completions of $Z$ with respect to these norms.
Then one can easily show that there exists a positive constant $C$ such that $\|\cdot\|_W\leq C\|\cdot\|_V$ on $Z$ so that there is a natural bounded linear operator $\mathcal{I}\colon V\rightarrow W$ extending the identity map on $Z$.
Since by~\cite{Anne89} -- see~\cite{ChadrakarGittinsHabibPeyerimhoff2025} for a corrected proof -- any closed and co-closed $p$-form vanishing along $\partial M$ must vanish identically on $M$, the operator $\mathcal{I}$ is actually injective.

Now let $\mathcal{K}\colon V\to V$ be the linear operator defined by
\[
\bigl(\mathcal{K}\omega,\omega'\bigr)_V=\bigl(\mathcal{I}\omega,
\mathcal{I}\omega'\bigr)_W
\]
for all $(\omega,\omega')\in V^2$.
By definition, the operator $\mathcal{K}$ is self-adjoint and positive-definite.
On the other hand, since from standard elliptic estimates both norms $\|\cdot\|_V$ and $\|\cdot\|_{H^2(M)}$ are equivalent on $Z$, the Rellich theorem ensures that $\mathcal{I}$ is compact and so is $\mathcal{K}$.
The spectral theorem for positive compact self-adjoint operators applies and yields the existence of a countable Hilbert or\-tho\-nor\-mal basis $(\omega_i)_{i\geq 1}$ of $V$ associated to a mo\-no\-to\-nous\-ly nonincreasing positive real sequence of eigenvalues of finite multiplicities $(\alpha_{i,p})_{i\geq1}$ going to $0$ such that $\mathcal{K}\omega_i=\alpha_{i,p}\omega_i$ for all $i\geq1$. Now fixing $i\geq 1$ and using the definition of $\mathcal{K}$ as well as the integration by part formula~(\ref{eq:laplacian}), it can be computed that, for every $\omega\in Z$,
\[%\label{eq:DeltaomegaDeltaomegaibis}
\alpha_{i,p}(\Delta\omega_i,\Delta\omega)_{L^2(M)} =(\mathcal{K}\omega_i,\omega)_V
=(\omega_i,\omega)_W=(\Delta\omega_i,\omega)_{L^2(M)}.
\]
At the same time, we also have by~\eqref{eq:partialintDeltacons} that
\[%\label{eq:DeltaomegaDeltaomegai}
(\Delta\omega_i,\Delta\omega)_{L^2(M)}= \bigl(\Delta^2\omega_i,\omega\bigr)_{L^2(M)}
\]
for every $\omega\in Z$ and so the previous equality now reads as
\[
\bigl(\Delta^2\omega_i-\Lambda_{i,p}\Delta \omega_i,
\omega\bigr)_{L^2(M)}=0,
\]
where we let \smash{$\Lambda_{i,p}=\frac{1}{\alpha_{i,p}}$}.
It follows that $\omega_i$ is a weak solution of the eigenvalue problem~(\ref{buckling_forms}) which, by ellipticity (see in Lemma~\ref{elliptic_problem}), is in fact smooth.
Thus the form $\omega_i$ becomes a smooth eigenform to problem~\eqref{buckling_forms} associated with the eigenvalue \smash{$\Lambda_{i,p}=\frac{1}{\alpha_{i,p}}$} which is of finite multiplicity, since $\alpha_{i,p}$ is.

Conversely, observe that if there exists a nontrivial solution $\omega$ to~\eqref{buckling_forms} for a certain $\Lambda\in\R$, then by~\eqref{eq:partialintDeltacons}, we have
\[
\bigl(\omega,\omega'\bigr)_V=\Lambda\bigl(\omega,\omega'\bigr)_W
\]
for every $\omega'\in Z$. Note that $\Lambda>0$, since otherwise $\Delta^2\omega=0$ which from~\eqref{eq:partialintDeltacons} implies that $\Delta\omega=0$ and then $\omega=0$ by~\cite{Anne89} since $\omega_{|_{\partial M}}=0$.
By definition of $\mathcal{K}$, we have $(\omega,\omega')_V=\Lambda(\mathcal{K}\omega,\omega')_V$ for all ${\omega'\in Z}$ and hence in $V$, therefore \smash{$\mathcal{K}\omega=\frac{1}{\Lambda}\omega$}. This shows that $\omega$ is an eigenform of $\mathcal{K}$ associated to the eigenvalue \smash{$\alpha=\frac{1}{\Lambda}$}.

Finally, given any eigenform $\omega$ associated to a positive eigenvalue $\Lambda$ of~\eqref{buckling_forms}, we have by formula~\eqref{eq:partialintDeltacons} that
\[
\Lambda\int_M\langle\Delta\omega,\omega\rangle {\rm d}\mu=\int_M |\Delta \omega|^2 {\rm d}\mu.
\]
Applying~\eqref{eq:laplacian} to the left-hand side of this equality ensures that
\[
\Lambda_{1,p}\leq \frac{\|\Delta\omega\|_{L^2(M)}^2}{
\|{\rm d}\omega\|_{L^2(M)}^2+\|\delta\omega\|_{L^2(M)}^2}
\]
for every such eigenform, with equality for $\omega$ associated to $\Lambda_{1,p}$.
Finally, if $\omega\in V$, one may write its decomposition in the Hilbert basis $(\omega_i)_{i\geq 1}$ so that
\begin{align*}
\|{\rm d}\omega\|_{L^2(M)}^2+\|\delta\omega\|_{L^2(M)}^2 &=(K\omega,\omega)_V
=\sum_{i\geq 1}\frac{1}{\Lambda_{i,p}}|(\omega,\omega_i)_V|^2
 \leq \frac{1}{\Lambda_{1,p}}\sum_{i\geq 1}|(\omega,\omega_i)_V|^2\\
 &=\frac{1}{\Lambda_{ 1 ,p}}\|\Delta\omega\|_{L^2(M)}^2.
\end{align*}
This prove the characterization~\eqref{eq:charb} since $Z$ is dense in $V$.
\end{proof}

\begin{Remark}\label{SymmetryBuck}
When $M$ is oriented, the Hodge $\star$ operator is an isometry commuting with the Laplacian and preserving the boundary conditions in the buckling problem so that $\Lambda_{i,p}=\Lambda_{i,n-p}$ for any $i\geq 1$ and $1\leq p\leq n$.
\end{Remark}

\subsection{The clamped plate problem}

The clamped plate eigenvalue problem on differential forms is
\begin{equation}\label{clamped_plate_forms}
\begin{aligned}
&\Delta^2 \omega=\Gamma \omega \ \quad&&\textrm{on } M,\\
&\omega =0 \ \quad&&\textrm{on }\partial M,&\\
&\nabla_\nu\omega_{|\partial M}=0 \ \quad &&\textrm{on }\partial M,
\end{aligned}
\end{equation}
for some real constant $\Gamma$. As previously, we immediately get the following existence result.
\begin{Theorem}\label{CP_Pb}
There exists a Hilbert basis of the space of $L^2$-integrable $p$-forms on $(M^n,g)$ consisting of eigenforms solutions of the problem~\eqref{clamped_plate_forms} associated to an unbounded and positive sequence of eigenvalues $(\Gamma_{i,p})_{i\geq 1}$. Moreover, each eigenspace has a finite multiplicity and the corresponding eigenforms are smooth. Finally, the first eigenvalue $\Gamma_{1,p}$ is characterized by
\begin{equation}\label{eq:charcp}
\Gamma_{1,p}=\inf\Biggl\{\frac{\|\Delta\omega\|_{L^2(M)}^2}{
\| \omega\|_{L^2(M)}^2},\, \omega\in\Omega^p(M)\setminus\{0\},\, \omega_{|_{\partial
M}}=0\textrm{ and }\nabla_\nu\omega_{|\partial M}=0\Biggr\}.
\end{equation}
Equality holds if and only if $\omega$ is an eigenform associated to the first eigenvalue.
\end{Theorem}

\begin{proof} It is enough to take $(\,,\,)_W$ to be the $L^2(M)$-scalar product on $\Omega^p(M)$ in the proof of Theorem~\ref{BucKPb} and then the proof goes the same. Note that if $\omega\in\Omega^p(M)$ is an eigenform associated to $\Gamma_{1,p}$, it follows from~(\ref{eq:partialintDeltacons}) that $\Gamma_{1,p}\geq 0$. Moreover, if $\Gamma_{1,p}=0$, then any associated eigenform $\omega$ has to be harmonic with $\omega_{|\partial M}=0$ and so $\omega=0$ by~\cite{Anne89}. In particular, $\Gamma_{1,p}>0$.
\end{proof}

\begin{Remark}
As for the buckling problem, when $M$ is oriented, the Hodge $\star$ operator preserves the clamped plate problem so that $\Gamma_{i,p}=\Gamma_{i,n-p}$ for any $i\geq 1$ and $1\leq p\leq n$.
\end{Remark}

\section{Eigenvalues of the buckling and clamped plate operators}\label{s:ebs}

\subsection{Eigenvalues for bounded Euclidean domains}%\label{ss:eigenvboundedEucldomains}

In this subsection, we completely describe the spectrum of the buckling and clamped plate problems for bounded domains in the Euclidean space. More precisely, we prove the following characterization.

\begin{Proposition}\label{BCP_Euclidean}
Let $(M^n,g)$ be a bounded domain in the Euclidean space $\mathbb{R}^n$.
Then the spectrum without multiplicities of the buckling problem on $p$-forms on $(M^n,g)$ coincides with the spectrum of the buckling problem on functions that is $\Lambda_{i,p}=\Lambda_{i,0}$ for all $i \geq 1$ and $p\in\{1,\dots,n\}$.
The same holds for the clamped plate problem that is $\Gamma_{i,p}=\Gamma_{i,0}$ for all $i \geq 1$ and $p\in\{1,\dots,n\}$.
The multiplicities of $\Lambda_{i,p}$ and $\Gamma_{i,p}$ are exactly $\bigl(\begin{smallmatrix}n\\ p\end{smallmatrix}\bigr)$ times those of $\Lambda_{i,0}$ and $\Gamma_{i,0}$ respectively.
\end{Proposition}

\begin{proof} First recall that on $\mathbb{R}^n$, there exists for each $p\in\{1,\dots,n\}$ a maximal number of parallel $p$-forms.
Fix $p\in\{1,\dots,n\}$ and denote by $\omega_0$ a nontrivial parallel $p$-form on $M$.
Then note that for any smooth function $f$ on $M$ with $f_{|\partial M}=0$ and \smash{$\frac{\partial f}{\partial \nu}_{|\partial M}=0$}, the $p$-form $\omega_f:=f \omega_0$ satisfies
\begin{equation}\label{fomega0BC}
\omega_{f_{|_{\partial M}}}=0\qquad\text{and}\qquad \nabla_\nu\omega_{f_{|_{\partial M}}}=0.
\end{equation}
On the other hand, since $\omega_0$ is parallel, we have ${\rm d}\omega_f={\rm d}f\wedge\omega_0$ and $\delta\omega_f=-{\rm d}f\lrcorner\omega_0$ and therefore $\Delta\omega_f=(\Delta f)\omega_0$. Applying twice this formula leads to
\begin{equation}\label{fomega0BL}
\Delta^2\omega_f=\Delta ((\Delta f)\omega_0 )=\bigl(\Delta^2 f\bigr)\omega_0.
\end{equation}
Now if we take $f_1$ (resp.\ $f_2$) to be an eigenfunction for the buckling (resp.\ clamped plate) problem~(\ref{buckling_functions}) (resp.\ (\ref{clamped_plate_functions})) associated with the eigenvalue $\Lambda$ (resp.\ $\Gamma$), we conclude combining~(\ref{fomega0BC}) and~(\ref{fomega0BL}) that $\omega_{f_1}$ (resp.\ $\omega_{f_2}$) is a $p$-eigenform for~(\ref{buckling_forms}) (resp.\ (\ref{clamped_plate_forms})) associated with the eigenvalue $\Lambda$ (resp.\ $\Gamma$).

Conversely, first note that if $f\in C^\infty(M)$ is a smooth nontrivial function such that $f=\<\omega,\omega_0\>$ where $\omega$ is a smooth $p$-form and $\omega_0$ is a smooth parallel $p$-form, we have
\begin{equation}\label{DeltaFunctionEuc}
\Delta f=\<\nabla^\ast\nabla\omega,\omega_0\>=\<\Delta\omega,\omega_0\>,
\end{equation}
where the last equality follows from the Bochner formula (see~(\ref{eq:weitzenboeckpforms}) below) and the fact that $M$ is Euclidean.
Note also that if $\omega$ satisfies the boundary condition~(\ref{Zcharacterization}) then $f$ and $\partial f/\partial\nu$ vanish on $\partial M$.
Now if $\omega_1$ and $\omega_2$ denote respectively $p$-eigenforms to the problems~(\ref{buckling_forms}) and~(\ref{clamped_plate_forms}) associated to the eigenvalues $\Lambda$ and $\Gamma$, then there exist two parallel $p$-forms $\omega^i_0$ on $\mathbb{R}^n$ such that $f_i:=\bigl\langle\omega_i,\omega^i_0\big\rangle$ are smooth nontrivial functions for $i=1,2$. Therefore, we easily deduce from~(\ref{DeltaFunctionEuc}) that these functions are smooth eigenfunctions of~(\ref{buckling_functions}) and~(\ref{clamped_plate_functions}) respectively associated to $\Lambda$ and $\Gamma$.
\end{proof}

\begin{Remarks}\quad
\begin{enumerate}\itemsep=0pt
\item Proposition~\ref{BCP_Euclidean} gives immediately the value of the first eigenvalues of the buckling and clamped plate problems on $p$-forms for the Euclidean ball of radius \smash{$R=\frac{1}{H_0}$}.
More precisely, it follows from~\cite[Section~1]{AL} that
\[
\Lambda_{1,p}=j_{\frac{n}{2},1}^2 H_0^2\qquad\text{and}\qquad\Gamma_{1,p}=k_{\frac{n}{2}-1,1}^4 H_0^4,
\]
where \smash{$j_{\frac{n}{2},1}$} is the first positive zero of the Bessel function \smash{$J_{\frac{n}{2}}$} and \smash{$k_{\frac{n}{2}-1,1}$} is the first positive zero of \smash{$J_{\frac{n}{2}-1}I_{\frac{n}{2}}+J_{\frac{n}{2}} I_{\frac{n}{2}-1}$}, with $I_\ell$ being the corresponding modified Bessel function of the~first kind.

\item If $(M^n,g)$ is a compact Riemannian manifold carrying a nontrivial parallel $p$-form for a~certain $p\in\{1,\dots,n\}$, we can mimic the first part of the proof of Proposition~\ref{BCP_Euclidean} to ensure that if ${\rm Spec}_{\Lambda,p}(M)$ and ${\rm Spec}_{\Gamma,p}(M)$ denote respectively the buckling and clamped plate spectra on $p$-forms for $p\in\{0,\dots,n\}$, then we have
\[
{\rm Spec}_{\Lambda,0}(M)\subset {\rm Spec}_{\Lambda,p}(M) \qquad\text{and}\qquad {\rm Spec}_{\Gamma,0}(M)\subset {\rm Spec}_{\Gamma,p}(M)
\]
\end{enumerate}
\end{Remarks}

\subsection{General estimates}

In this subsection, we prove general estimates between the first eigenvalues of the buckling and the clamped plate problems and other classical problems. These estimates hold on $(M^n,g)$ without any extra assumption, except being a compact Riemannian manifold with smooth boundary.

The first result of this part establishes a direct link between the buckling and the clamped plate problems.
\begin{Theorem}\label{BuckCP}
On an $n$-dimensional compact Riemannian manifold $(M^n,g)$ with smooth boundary, we have $\Gamma_{1,p}<\Lambda_{1,p}^2$ for all $p\in\{0,\dots,n\}$.
\end{Theorem}

Before giving the proof of this result, we first note that if $\omega$ is a $p$-form satisfying $\omega_{|_{\partial M}}=\nabla_\nu \omega_{|_{\partial M}}=0$, then it follows from the identity $\Delta={\rm d}\delta+\delta {\rm d}$ and from an integration by parts~that
\[%\label{normDeltaomega}
\|\Delta \omega \|^2_{L^2(M)} = \|\delta {\rm d}\omega \|^2_{L^2( M)}+\|{\rm d}\delta\omega \|^2_{L^2(M)}
\]
and from~(\ref{eq:laplacian}) that
\begin{equation}\label{Deltaomega_scal_omega}
(\Delta \omega,\omega)_{L^2( M)} = \|{\rm d}\omega \|^2_{L^2( M)}+\|\delta\omega \|^2_{L^2(M)}.
\end{equation}

\begin{proof} Let $\omega$ be a smooth nonzero $p$-form on $M$ with $\omega_{|_{\partial M}}=0$ and $\nabla_\nu\omega_{|\partial M}=0$. By~\eqref{Deltaomega_scal_omega} and the Cauchy--Schwarz inequality, we get
\[
\|{\rm d}\omega \|^2_{L^2( M)}+\|\delta\omega \|^2_{L^2( M)}\leq\| \Delta \omega\|_{L^2( M)} \| \omega\|_{L^2( M)}
\]
and this implies that
\[
\frac{\| \Delta \omega\|_{L^2( M)}}{\|\omega\|_{L^2(M)}}\leq\frac{\| \Delta \omega\|_{L^2( M)}^2}{\|{\rm d}\omega \|^2_{L^2( M)}+\|\delta\omega \|^2_{L^2( M)}}.
\]
Taking $\omega$ to be an eigenform of the buckling problem associated to $\Lambda_{1,p}$ gives the desired inequality. Moreover, if equality holds, then for any eigenform $\omega$ of the buckling problem associated to~$\Lambda_{1,p}$, the $p$-form $\Delta\omega$ is pointwise proportional to $\omega$, therefore $\Delta\omega=\Lambda_{1,p}\omega$ on $M$ and thus~$\omega$ is an eigenform of the Dirichlet problem associated to the eigenvalue $\Lambda_{1,p}$.
But because $\omega$ satisfies $\nabla_\nu\omega_{|\partial M}=0$ along $\partial M$, the unique continuation theorem for elliptic second-order linear operators (see, e.g.,~\cite[Theorem~1.4]{Salolectnotes2014} and~\cite[Chapter~VIII]{HoermanderlinearI}) implies that $\omega=0$ on $M$. This shows that the inequality must actually be strict.
\end{proof}

Now we consider $\lambda_{1,p}$ the smallest eigenvalue of the Hodge Laplacian of $M$ with Dirichlet boundary condition.
It is well-known that there exists a Hilbert basis of the space of $L^2$-integrable $p$-forms on $(M^n,g)$ consisting of smooth eigenforms solutions of the problem
\begin{equation}\label{DirichletBC}
\begin{aligned}
&\Delta \omega =\lambda \omega \ \quad&&\textrm{on } M,\\
&\omega=0\ \quad &&\textrm{on }\partial M
\end{aligned}
\end{equation}
associated to an unbounded and positive sequence of eigenvalues $(\lambda_{i,p})_{i\geq 1}$. Moreover the first eigenvalue is characterized by
\[
\lambda_{1,p}=\inf\Biggl\{\frac{\|{\rm d}\omega\|_{L^2(M)}^2+\|\delta\omega\|_{L^2(M)}^2}{\|\omega\|_{L^2(M)}^2},\, \omega\in\Omega^p(M)\setminus\{0\},\, \omega_{|_{\partial M}}=0\Biggr\}
\]
and has to be positive because of~\cite{Anne89}.

\begin{Remarks}\quad%\label{RemarksDHL}
\begin{enumerate}\itemsep=0pt
\item The Hodge $\star$ operator preserves the Dirichlet problem so that $\lambda_{i,p}=\lambda_{i,n-p}$ for any $i\geq 1$ and $1\leq p\leq n$.
\item If $(M^n,g)$ is a bounded domain in the Euclidean space, one can show, reasoning as in the~proof of Proposition~\ref{BCP_Euclidean}, that the Dirichlet eigenvalues on $p$-forms do not depend on $p$ and correspond to the Laplacian Dirichlet eigenvalues on functions.
\end{enumerate}
\end{Remarks}

We are now in position to give explicit estimates between the first eigenvalues of the three previous problems, restricting our attention to \smash{$1\leq p\leq[\frac{n}{2}]$} by Hodge symmetry of these spectra.
Namely, we get the following.

\begin{Theorem}\label{BCP_Dirichlet}
On an $n$-dimensional compact Riemannian manifold $(M^n,g)$ with smooth boundary, the following inequalities hold:
\begin{enumerate}\itemsep=0pt
\item[$(1)$] $\Lambda_{1,p} \lambda_{1,p}< \Gamma_{1,p}$ for $0\leq p\leq n$,
\item[$(2)$] $\inf(\lambda_{1,p+1}, \lambda_{1,p-1}) \le \Lambda_{1,p}$ for $1\leq p\leq n$,
\item[$(3)$] $\lambda_{1,1}\leq \Lambda_{1,0}$.
\end{enumerate}
\end{Theorem}

\begin{Remarks}\quad
\begin{enumerate}\itemsep=0pt
\item The inequalities of Theorems~\ref{BuckCP} and~\ref{BCP_Dirichlet} give the analogue of some well-known results in the case of functions (see~\cite{AL} for example).
\item In general, it is not clear whether equality can occur in the second and third inequalities of Theorem~\ref{BCP_Dirichlet}.
We will see that this cannot occur for bounded Euclidean domains.
\end{enumerate}
\end{Remarks}

\begin{proof} To prove the first inequality, we consider $\omega$ an eigenform of the clamped plate problem~\eqref{clamped_plate_forms} associated with $\Gamma_{1,p}$. Then the form $\omega$ can be considered as a test-form for the first eigenvalue of the buckling problem as well as for the first eigenvalue of the Dirichlet problem on differential forms. Therefore, we first get
\[
\Lambda_{1,p}\bigl(\|{\rm d}\omega\|^2_{L^2(M)} +\|\delta \omega\|^2_{L^2(M)}\bigr) \le \|\Delta \omega\|^2_{L^2(M)},
\]
which can be rewritten as
\[
\Lambda_{1,p} \Biggl(\frac{\|{\rm d}\omega\|^2_{L^2(M)} +\|\delta \omega\|^2_{L^2(M)}}{\| \omega\|^2_{L^2(M)}}\Biggr)\le \frac{\|\Delta \omega\|^2_{L^2(M)}}{\| \omega\|^2_{L^2(M)}},
\]
from which the estimate follows. If the inequality is an equality, then any eigenform $\omega$ of the clamped plate problem must be an eigenform of the Dirichlet problem for the Hodge Laplacian on $M$.
As above, the unique continuation theorem implies that $\omega=0$.
Therefore, the inequality is strict.

For the second estimate, we consider $\omega$ an eigenform of the buckling eigenvalue problem~(\ref{buckling_forms}) associated with the first eigenvalue $\Lambda_{1,p}$ for \smash{$1\leq p\leq [\frac{n}{2}]$}. We first notice that the two differential forms $\omega_1={\rm d}\omega\in\Omega^{p+1}(M)$ and $\omega_2=\delta\omega\in\Omega^{p-1}(M)$ cannot be both trivial otherwise $\omega$ would be a harmonic $p$-form which vanishes on the boundary which is impossible by~\cite{Anne89}. Moreover, it~follows from~(\ref{Zcharacterization}) that
\[
\iota^*\omega_1={\rm d}\iota^*\omega=0,\qquad \nu\lrcorner\omega_1=\nu\lrcorner {\rm d}\omega=0
\]
and
\[
\iota^*\omega_2=\iota^*(\delta\omega)=0,\qquad\nu\lrcorner\omega_2=-\delta^{\partial M}(\nu\lrcorner\omega)=0,
\]
which imply that both $\omega_1$ and $\omega_2$ vanish along the boundary. So one can respectively take these forms as test-forms in~(\ref{DirichletBC}) leading to
\[
\lambda_{1,p+1}\| {\rm d} \omega\|^2_{L^2(M)} \le \|\delta {\rm d} \omega\|^2_{L^2(M)} \qquad \text{and} \qquad \lambda_{1,p-1}\| \delta \omega\|^2_{L^2(M)} \le \|{\rm d}\delta \omega\|^2_{L^2(M)}.
\]
By adding these two inequalities, we get
\[
\inf(\lambda_{1,p+1}, \lambda_{1,p-1}) \bigl(\| {\rm d} \omega\|^2_{L^2(M)} + \| \delta \omega\|^2_{L^2(M)} \bigr) \le \|\delta {\rm d} \omega\|^2_{L^2(M)}+\|{\rm d}\delta \omega\|^2_{L^2(M)}
\]
and the variational characterization~\eqref{eq:charb} allows to conclude.

For the last one, take a smooth eigenfunction $f$ for the buckling problem~(\ref{buckling_functions}) on functions associated to $\Lambda_{1,0}$ and let $\omega={\rm d}f\in\Omega^1(M)$. Then as above we observe that $\omega$ vanishes on $\partial M$ so that it can be used as a test-form in the variational characterization of $\lambda_{1,1}$ leading to the inequality
\[
\lambda_{1,1}\leq \frac{\int_M\bigl(\Delta f)^2{\rm d}\mu}{\int_M |{\rm d}f|^2{\rm d}\mu}=\Lambda_{1,0}.
\]
This proves the third inequality in the broad sense.
\end{proof}

A direct consequence of Theorem~\ref{BuckCP} and the first inequality in Theorem~\ref{BCP_Dirichlet} is the following estimate.
\begin{Corollary}\label{DBCP}
Let $(M^n,g)$ be a compact Riemannian manifold with boundary, then
\begin{equation}\label{DirichletBounds}
\lambda_{1,p}<\sqrt{\Gamma_{1,p}}<\Lambda_{1,p}
\end{equation}
for all $0\leq p\leq n$.
\end{Corollary}

Since for bounded domains in Euclidean space, the eigenvalues $\lambda_{1,p}$ and $\Lambda_{1,p}$ do not depend on $p$, it follows from the previous corollary with $p=0$ that equality cannot occur in the second and third inequalities of Theorem~\ref{BCP_Dirichlet}.

\begin{Remark} In~\cite{elchamiginouxhabib19}, the Robin problem is defined as
\[
\begin{aligned}
&\Delta \omega =\lambda \omega \ \quad&&\textrm{on } M, \\
&\nu\lrcorner\omega =0\ \quad&&\textrm{on }\partial M, \\
&\nu\lrcorner {\rm d}\omega = \tau \iota^*\omega \ \quad && \textrm{on }\partial M
\end{aligned}
\]
and it is proven that the first eigenvalue satisfies $\lambda_{1,p}(\tau)\leq \lambda_{1,p}$ for any parameter $\tau>0$.
Hence, it follows from~\eqref{DirichletBounds} that \smash{$\lambda_{1,p}(\tau)<\sqrt{\Gamma_{1,p}}<\Lambda_{1,p}$}.
\end{Remark}

\begin{Remark}
It is in fact not difficult to show that $\lambda_{j,p}\leq\Lambda_{j,p}$ for all $j\geq 1$ and $p=0,\dots,n$, where the inequality is strict for $j=1$ by Corollary~\ref{DBCP}.
Indeed, it can be checked that the eigenvalue problem~(\ref{DirichletBC}) is equivalent to the problem
\[
\begin{aligned}
&\Delta^2 \omega =\lambda \Delta\omega \ \quad&&\textrm{on } M, \\
&\omega =0 \ \quad&&\textrm{on }\partial M, \\
&\Delta\omega = 0 \ \quad && \textrm{on }\partial M
\end{aligned}
\]
in such a way that the Dirichlet eigenvalues are determined by applying the min-max formula to the functional
\[
\frac{\int_M|\Delta \omega|^2{\rm d}\mu}{\int_M|{\rm d}\omega|^2+|\delta\omega|^2{\rm d}\mu}
\]
 on subspaces of $p$-forms which are in $\bigl(H^2\cap H_0^1\bigr)(M)$.
Namely, any critical point $\omega\in\bigl(H^2\cap H_0^1\bigr)(M)$ of that functional associated with a critical value $\ell$ must satisfy
\[
\int_M\bigl\langle\Delta\omega,\Delta\omega'\bigr\rangle {\rm d}\mu=\ell\cdot\int_M\bigl(\bigl\langle {\rm d}\omega,{\rm d}\omega'\bigr\rangle+\bigl\langle \delta\omega,\delta\omega'\bigr\rangle\bigr){\rm d}\mu
\]
for all $\omega'\in\bigl(H^2\cap H_0^1\bigr)(M)$, that is,
\[
\int_M\bigl\langle\Delta^2\omega,\omega'\bigr\rangle {\rm d}\mu+\int_{\partial M}\bigl(\bigl\langle\iota^*\Delta\omega,\nu\lrcorner {\rm d}\omega'\bigr\rangle-\bigl\langle\nu\lrcorner\Delta\omega,\iota^*\delta\omega'\bigr\rangle\bigr){\rm d}\mu=\ell\cdot\biggl(\int_M\Delta\omega,\omega'\bigr\rangle\biggr),
\]
which reduces to
\[
\int_M\bigl\langle\Delta^2\omega-\ell\Delta\omega,\omega'\bigr\rangle {\rm d}\mu=\int_{\partial M}\bigl(\bigl\langle\nu\lrcorner\Delta\omega,\iota^*\delta\omega'\bigr\rangle-\bigl\langle\iota^*\Delta\omega,\nu\lrcorner {\rm d}\omega'\bigr\rangle\bigr){\rm d}\mu.
\]
Here we have used that $\omega'$ vanishes along $\partial M$.
It can be deduced in the usual way that ${\Delta^2\omega-\ell\Delta\omega=0}$ on $M$ as well as $\nu\lrcorner\Delta\omega=0$ and $\iota^*\Delta\omega=0$ along $\partial M$, using the fact that, for $\omega'\in H_0^1(M)$, we have $\iota^*\delta\omega'=-\nu\lrcorner\nabla_\nu\omega'$ and $\nu\lrcorner {\rm d}\omega'=\iota^*\nabla_\nu\omega'$, and both the tangential and normal components of $\nabla_\nu\omega'$ can be prescribed independently.
Thus $\Delta\omega=0$ along $\partial M$, which yields the above equivalent problem.

Comparing now the above characterization to that of $\Lambda_{j,p}$, where the same quotient is considered on the smaller space $H^2_0(M)$, allows to conclude.
Mind that it cannot be deduced for~${j\geq2}$ whether the equality can be attained or not.
\end{Remark}

Now for $p\in\{0,\dots,n\}$, we consider the Hodge Laplacian with respect to the absolute boundary condition which satisfies the eigenvalue problem
\begin{equation}\label{eq:Neumannpblonpforms}
\begin{aligned}
&\Delta \omega =\mu \omega\ \quad &&\textrm{on } M, \\
&\nu\lrcorner\omega =0\ \quad&&\textrm{on }\partial M, \\
&\nu\lrcorner\,{\rm d}\omega = 0\ \quad &&\textrm{on }\partial M.
\end{aligned}
\end{equation}

We will denote by $\mu_{1,p}$ the first {\it positive} eigenvalue which is given by
\begin{equation}\label{VC_Absolute}
\mu_{1,p}=\inf_{\omega\in\Omega^p(M)\setminus\{0\}}\Biggl\{\frac{\|{\rm d}\omega\|_{L^2(M)}^2+\|\delta\omega\|_{L^2(M)}^2}{\|\omega\|_{L^2(M)}^2},\ \nu\lrcorner\omega=0\textrm{ and }\omega\in H_A^p(M)^\perp\Biggr\}.
\end{equation}

The eigenspace associated with the zero eigenvalue, if non empty, corresponds to the absolute de Rham cohomology group in degree $p$ defined by
\[
H_A^p(M)= \{\omega\in \Omega^p(M)\mid {\rm d}\omega=\delta\omega=0\ \text{on}\ M\text{ and } \nu\lrcorner\omega=0 \}.
\]

With these notations, it is clear that $\mu_{1,0}$ corresponds to the first nonzero eigenvalue of the Laplace operator under the Neumann boundary condition. Note also that the dual eigenvalue problem is the relative one whose zero eigenvalue reflects the relative de Rham cohomology group
\[
 H_R^p(M)= \{\omega\in \Omega^p(M)\mid {\rm d}\omega=\delta\omega=0\ \text{on}\ M\text{ and } \iota^*\omega=0 \}
\]
and is related to the absolute boundary condition by the Hodge star operator which induces an~isomorphism \smash{$H_A^p(M)\simeq H_R^{n-p}(M)$} for all $p=0,\dots,n$.
The Hodge star operator maps~(\ref{eq:Neumannpblonpforms}) onto the boundary problem
\[
\begin{aligned}
&\Delta \omega =\kappa \omega\ \quad &&\textrm{on } M,\\
&\iota^*\omega =0 \ \quad &&\textrm{on }\partial M,\\
&\iota^*(\delta\omega) = 0\ \quad &&\textrm{on }\partial M,
\end{aligned}
\]
where $\kappa$ is the corresponding eigenvalue.
As usual, up to the possible eigenvalue $0$, the (mo\-no\-to\-nous\-ly nondecreasing) sequence of real eigenvalues is denoted by $(\kappa_{i,p})_{i\geq1}$, where $\kappa_{1,p}$ is the smallest positive one.
It is a straightforward consequence of that correspondence via the Hodge star operator that, for all $i\geq1$ and $p\in\{0,\dots,n\}$, the identity $\mu_{i,p}=\kappa_{i,n-p}$ holds.
Moreover, $\mu_{1,n}=\kappa_{1,0}>0$ holds because of \smash{$H^0_R(M)=\{0\}$}. By definition, the identity $\kappa_{1,0}=\lambda_{1,0}$ also~holds.

\begin{Remark}\label{RemarkPolya}
If $H_A^p(M)=\{0\}$, an eigenform for the Dirichlet problem~(\ref{DirichletBC}) can be taken as a~test-form in the variational characterization of $\mu_{1,p}$ leading to the estimate $\mu_{1,p}\leq\lambda_{1,p}$ for all $p=1,\dots,n$.
Note that this inequality does not hold in general for $p=0$.
However, it is a~well-known result of P\'olya~\cite{Polya1952} (see also~\cite{Friedlander91}) that $\mu_{1,0}<\lambda_{1,0}$ for bounded Euclidean domains.
On the other hand, for $p=0,\dots,[\frac{n}{2}]$, Guerini and Savo proved in~\cite[Theorem 2.6\,(b))]{GueriniSavo2004} that~${\mu_{1,p}\leq\kappa_{1,p}}$ for a~convex bounded domain $M$ in the Euclidean space.
Even more, if $M$ is strictly convex, this inequality is strict for \smash{$0\leq p<[\frac{n}{2}]$} (see~\cite[Theorem 2.6\,(c))]{GueriniSavo2004}).
Therefore, for a~bounded convex Euclidean domain with \smash{$H^{p}_R(M)=\{0\}$} for some \smash{$p=0,\dots,[\frac{n}{2}]$}, we have
\[
\mu_{1,p}\leq \kappa_{1,p}=\mu_{1,n-p}\leq\lambda_{1,n-p}=\lambda_{1,p}.
\]
This is the case, for example, for strictly convex bounded domains in $\mathbb{R}^n$ and in addition the first inequality is strict for \smash{$0\leq p<[\frac{n}{2}]$} by the above discussion.
\end{Remark}

It turns out that the first nonzero eigenvalue of the Hodge Laplacian under the absolute boundary condition can be related with the first eigenvalue of the buckling problem, as is stated below.

\begin{Theorem}\label{BuckAbs}
On an $n$-dimensional oriented compact Riemannian manifold $(M^n,g)$ with smooth boundary, we have
\begin{equation}\label{BuckAbs1}
\max(\mu_{1,p},\mu_{1,n-p})\leq\Lambda_{1,p}
\end{equation}
for all $p=0,\dots,n$. Moreover, equality occurs for some $p$ if and only if there exist $k\in\{p,n-p\}$, $\omega\in\Omega^k(M)$ as well as $\omega_0\in H_A^k(M)\setminus\{0\}$ satisfying the overdetermined boundary value problem%
\begin{equation}\label{ODP}
\begin{aligned}
&\Delta \omega =\mu_{1,k} \omega \ \quad&&\textrm{on } M, \\
&\omega = \omega_0 \ \quad&& \textrm{on }\partial M,\\
&\iota^*(\delta\omega) =0,\ \nu\lrcorner {\rm d}\omega =0 \ \quad&& \textrm{on }\partial M.
\end{aligned}
\end{equation}
\end{Theorem}

\begin{proof} Let $\omega_1\in\Omega^p(M)$ be an eigenform of the buckling problem associated to the eigenvalue~$\Lambda_{1,p}$.
We denote by $\omega_0$ the $L^2(M)$-orthogonal projection of $\omega_1$ onto \smash{$H_A^p(M)$} and we set $\omega:=\omega_1-\omega_0\in \Omega^p(M)$.
Note that $\omega_0$ may vanish, which is the case if and only if $\nu\lrcorner\,{\rm d}\Delta\omega_1$ is $L^2(\partial M)$-orthogonal to $\iota^*H_A^p(M)$ by~\eqref{eq:partialintegDelta2omega}.
Then ${\rm d}\omega={\rm d}\omega_1$ and $\delta\omega=\delta\omega_1$ on $M$, and in particular $\Delta\omega=\Delta\omega_1$, as well as $\nu\lrcorner\,\omega=0$ and $\nu\lrcorner\,{\rm d}\omega=0$ along $\partial M$.
For this form $\omega$, using the first equality in~(\ref{eq:laplacian}) and the Cauchy--Schwarz inequality leads to
\begin{gather*}
\bigl(\|{\rm d}\omega\|_{L^2(M)}^2+\|\delta\omega\|_{L^2(M)}^2\bigr)^2 = (\Delta\omega,\omega )_{L^2(M)}^2 \leq \|\Delta\omega\|_{L^2(M)}^2 \|\omega\|_{L^2(M)}^2,
\end{gather*}
so that, together with the fact that $\omega\neq0$ (otherwise, $\Lambda_{1,p}=0$),
\[
\frac{\|{\rm d}\omega\|_{L^2(M)}^2+\|\delta\omega\|_{L^2(M)}^2}{\|\omega\|_{L^2(M)}^2}\leq\frac{\|\Delta\omega\|_{L^2(M)}^2}{\|{\rm d}\omega\|_{L^2(M)}^2+\|\delta\omega\|_{L^2(M)}^2}.
\]
Note that \smash{$\|{\rm d}\omega\|_{L^2(M)}^2+\|\delta\omega\|_{L^2(M)}^2=\|{\rm d}\omega_1\|_{L^2(M)}^2+\|\delta\omega_1\|_{L^2(M)}^2>0$}.
Moreover, since $\Delta\omega=\Delta\omega_1$ and $\omega_1$ is an eigenform associated with $\Lambda_{1,p}$, we deduce that
\[
\frac{\|{\rm d}\omega\|_{L^2(M)}^2+\|\delta\omega\|_{L^2(M)}^2}{\|\omega\|_{L^2(M)}^2}\leq\frac{\|\Delta\omega_1\|_{L^2(M)}^2}{\|{\rm d}\omega_1\|_{L^2(M)}^2+\|\delta\omega_1\|_{L^2(M)}^2}=\Lambda_{1,p}.
\]
But $\nu\lrcorner\omega=0$ and \smash{$(\omega,\omega')_{L^2(M)}=0$} for all $\omega'\in H_A^p(M)$ so that the $p$-form $\omega$ can be taken as a~test-form in~(\ref{VC_Absolute}) and therefore $\displaystyle\mu_{1,p}\leq \Lambda_{1,p}$.
The same argument ensures that $\mu_{1,n-p}\leq\Lambda_{1,n-p}$ and the result follows from $\Lambda_{1,n-p}=\Lambda_{1,p}$ as was noticed in Remark~\ref{SymmetryBuck}.

If equality occurs, then it can be assumed without loss of generality that $\mu_{1,p}=\Lambda_{1,p}$ up to replacing $p$ by $n-p$.
Then the $p$-form $\omega$ is a $\mu_{1,p}$-Laplace-eigenform with respect to the absolute boundary condition, so that it satisfies the problem~(\ref{eq:Neumannpblonpforms}).
On the other hand, it is clear that, since $\omega=\omega_1-\omega_0$
%\footnote{Would $\omega_1-t\omega_0$ be better-suited?}
with $\omega_1$ a buckling-eigenform and $\omega_0\in H^p_A(M)$, we have $\iota^*\omega=-\iota^*\omega_0$ and ${\iota^*(\delta\omega)=0}$, therefore $\omega$ solves~\eqref{ODP} where $-\omega_0$ takes the role of $\omega_0$.
Note that necessarily the $\omega_0$ from~\eqref{ODP} cannot vanish, otherwise $\omega=0$ by the unique continuation property for second-order elliptic systems, which would be a contradiction.
This gives the first part of the equality case.
Conversely, assume the existence of $\omega\in\Omega^p(M)$ satisfying~(\ref{ODP}) with $\omega_0\in H^p_A(M)\setminus\{0\}$.
Let $\omega_1:=\omega-\omega_0\in\Omega^p(M)$.
Notice that $\omega_1\neq0$ (otherwise $\mu_{1,p}=0$) and that $\omega$ is $L^2(M)$-orthogonal to $H_A^p(M)$ by integrating $\Delta\omega=\mu_{1,p}\omega$ against any $\omega_0'\in H_A^p(M)$ on $M$ and using~\eqref{eq:laplacianG1}.
Then it is not difficult to check using~(\ref{Zcharacterization}) that
\[
\omega_{1}{}_{|_{\partial M}}=0\qquad\text{and}\qquad\nabla_\nu\omega_{1}{}_{|_{\partial M}}=0,
\]
so that $\omega_1$ can be taken as a test-form in the variational characterization of $\Lambda_{1,p}$.
On the other hand, since $\omega_0$ is a harmonic form, we have ${\rm d}\omega_1={\rm d}\omega$ and $\delta\omega_1=\delta\omega$ and then
\[
\Lambda_{1,p}\leq\frac{\|\Delta\omega_1\|^2_{L^2(M)}}{\|{\rm d}\omega_1\|^2_{L^2(M)}+\|\delta\omega_1\|^2_{L^2(M)}}=\frac{\|\Delta\omega\|^2_{L^2(M)}}{\|{\rm d}\omega\|^2_{L^2(M)}+\|\delta\omega\|^2_{L^2(M)}}=\mu_{1,p}.
\]
This implies that $\mu_{1,p}=\Lambda_{1,p}$.
Note that~\eqref{ODP} is invariant under the Hodge star operator: if~${\omega\in\Omega^p(M)}$ solves~\eqref{ODP}, then so does $\star\omega\in\Omega^{n-p}(M)$ for the same $\mu=\mu_{1,p}$.
This shows that, if $\mu_{1,p}=\Lambda_{1,p}$, then actually $\mu_{1,p}$ must also be a Laplace-eigenvalue with absolute boundary conditions on $n-p$-forms, however it does not have to be the smallest one.
This concludes the~proof.
\end{proof}

Note that problem~(\ref{ODP}) would not be elliptic without the boundary condition $\nu\lrcorner\,{\rm d}\omega=0$.

\begin{Remark}
In other words, Theorem~\ref{BuckAbs} states that
\[%\label{BuckAbs2}
\max (\mu_{1,p},\kappa_{1,p} )\leq\Lambda_{1,p}
\]
for all $p=0,\dots,n$. Moreover, equality occurs if and only if there exist either $\omega\in\Omega^p(M)$ and \smash{$\omega_0\in H^p_A(M)\setminus\{0\}$} satisfying~(\ref{ODP}), or \smash{$\widetilde{\omega}\in\Omega^p(M)$} and \smash{$\widetilde{\omega}_0\in H^p_R(M)\setminus\{0\}$} such that
\begin{equation}\label{ODP2}
\begin{aligned}
&\Delta \widetilde{\omega} =\kappa_{1,p} \widetilde{\omega} \ \quad &&\textrm{on } M, \\
&\widetilde{\omega} = \widetilde{\omega}_0 \ \quad && \textrm{on }\partial M, \\
&\iota^*(\delta\widetilde{\omega}) =0,\ \nu\lrcorner {\rm d}\widetilde{\omega} =0\ \quad && \textrm{on }\partial M.
\end{aligned}
\end{equation}
For $p=0$ or $p=n$,~\eqref{ODP2}
% becomes trivial (the only solution is an identically vanishing function) as mentioned in Theorem~\ref{BuckAbs}.
and~\eqref{ODP} are both equivalent to the existence of a nonzero function~$f$ satisfying $\Delta f=\mu_{1,0}f$ on $M$ with $f_{|_{\partial M}}=1$ (up to rescaling $f$) as well as $\partial_\nu f=0$ along $\partial M$.
The existence of a nontrivial such solution is related to the so-called Schiffer conjecture, see comments after Remarks~\ref{r:Schifferreformulated} below.
\end{Remark}

As a direct consequence of Theorem~\ref{BuckAbs}, we get the following.
\begin{Corollary}
Let $(M^n,g)$ be an $n$-dimensional compact Riemannian manifold with boundary, then the following hold:
\begin{enumerate}\itemsep=0pt\samepage
\item[$(1)$] if $H^p_A(M)=\{0\}$ then $\mu_{1,p}<\Lambda_{1,p}$ for $p=0,\dots,n$;
\item[$(2)$] if $H^p_A(M)=\{0\}$ and $H^p_R(M)=\{0\}$ then inequality~\eqref{BuckAbs1} is strict for $p=0,\dots,n$.
\end{enumerate}
\end{Corollary}

Note that, when \smash{$H^p_A(M)=\{0\}$}, the inequality $\mu_{1,p}<\Lambda_{1,p}$ is a straightforward consequence of Corollary~\ref{DBCP} because of $\max(\mu_{1,p},\mu_{1,n-p})\leq \lambda_{1,p}$ by~\cite[Proposition~5.4]{elchamiginouxhabib19}.

\begin{Remarks}\label{r:Schifferreformulated}
\noindent\begin{enumerate}\itemsep=0pt
\item Since \smash{$H^n_A(M)=\{0\}$}, it follows from Theorem~\ref{BuckAbs} that $\mu_{1,n}<\Lambda_{1,n}$, which directly follows from the inequality $\lambda_{1,0}<\Lambda_{1,0}$.
\item Let $\omega_i\in\Omega^p(M)$, $i=1,2$, be two $p$-forms on $M$ such that $\omega_{1|\partial M}=\omega_{2|\partial M}$. Then it follows directly from~(\ref{eq:iotastarnablanuomega}) and~(\ref{eq:nuintnablanuomega}) that
\[
\bigl(\iota^*(\delta\omega_1)=\iota^*(\delta\omega_2)\quad\text{and}\quad\nu\lrcorner {\rm d}\omega_1=\nu\lrcorner {\rm d}\omega_2\bigr)\qquad\Longleftrightarrow\qquad\nabla_\nu\omega_{1|\partial M}=\nabla_\nu\omega_{2|\partial M}.
\]
Using this characterization, it is not difficult to show that $\omega\in\Omega^p(M)$ and \smash{$\omega_0\in H^p_A(M)$} satisfy~(\ref{ODP}) if and only if
\[
\begin{aligned}
&\Delta \omega =\mu_{1,p} \omega \ \quad &&\textrm{on } M, \\
&\omega = \omega_0\ \quad && \textrm{on }\partial M, \\
&\nabla_\nu\omega_{|\partial M} = \nabla_\nu\omega_{0|\partial M}\ \quad &&\textrm{on }\partial M.
\end{aligned}
\]
Obviously, the same holds for the eigenvalue boundary problem~(\ref{ODP2}) with $\widetilde{\omega}\in\Omega^p(M)$ and \smash{$\widetilde{\omega}_0\in H^p_R(M)$}.
\end{enumerate}
\end{Remarks}

From the estimates~\eqref{BuckAbs1} and~\eqref{DirichletBounds}, we have, for $p=0$, both $\mu_{1,0}\leq\Lambda_{1,0}$ and $\lambda_{1,0}<\Lambda_{1,0}$.
% The second inequality is a consequence of .
The first inequality was first proved in~\cite{IliasShouman2019} and states more precisely that $\mu_{1,0}<\Lambda_{1,0}$ on any compact Riemannian manifold with boundary. However the arguments given in the proof of~\cite[Lem\-ma~3.1]{IliasShouman2019}, which establishes that the inequality is strict, are not clear.
Indeed, if $M$ is assumed to be connected, we have $H^0_A(M)\simeq \mathbb{R}$ and so, as indicated in Theorem~\ref{BuckAbs}, the equality $\mu_{1,0}=\Lambda_{1,0}$ ensures the existence of a smooth function $f$ on $M$ satisfying the overdetermined problem
\[
\begin{aligned}
&\Delta f =\mu_{1,0} f\ \quad  &&\textrm{on } M,\\
&f = 1\ \quad && \textrm{on }\partial M,\\
&\frac{\partial f}{\partial \nu} = 0 \ \quad &&\textrm{on }\partial M,
\end{aligned}
\]
which is a particular case of the so-called Schiffer conjecture (see, for example,~\cite{ProvenzanoSavo2023} and the references therein).
However, the fact that the inequality is strict can be proved at least for bounded Euclidean domains.
\begin{Corollary}%\label{c:mu10<Lambda10bddEucldomains}
Let $(M^n,g)$ be a bounded domain in the Euclidean space. Then $\mu_{1,0}<\Lambda_{1,0}$.
\end{Corollary}

\begin{proof} As recalled in Remark~\ref{RemarkPolya}, we have $\mu_{1,0}<\lambda_{1,0}$ for any bounded domain in Euclidean space and so the result follows from Corollary~\ref{DBCP}.
\end{proof}

On a bounded convex Euclidean domain with \smash{$H^{p}_R(M)=\{0\}$}, since the spectrum of the Dirichlet and the buckling problems do not depend on $p$, we deduce from Remark~\ref{RemarkPolya} that
\[
\mu_{1,p}\leq\lambda_{1,p}=\lambda_{1,0}<\lambda_{2,0}\leq\Lambda_{1,0}=\Lambda_{1,p}
\]
for \smash{$p=0,\dots,[\frac{n}{2}]$}.
The last inequality is due to Payne~\cite{Payne1955} (with a gap in the proof filled by Friedlander~\cite{Friedlander2004}).

\subsection{Manifolds with Weitzenb\"ock curvature operator bounded from below}

In this subsection, we give a lower bound for the first eigenvalue $\lambda_{1,p}$ of the Dirichlet Hodge Laplacian under the condition that the Weitzenb\"ock curvature operator of the Riemannian manifold $(M^n,g)$ is bounded from below by a positive constant $\gamma>0$. Combined with some results of the former part, we obtain estimates for the first eigenvalues of the buckling and clamped plate problems in this context. These last results generalize previous results by Chen, Cheng, Wang and Xia~\cite{ChenChengWangXia2012}.

For this, we first recall that the Weitzenb\"ock formula for $p$-forms states that
\begin{equation}\label{eq:weitzenboeckpforms}
\Delta\omega=\nabla^*\nabla\omega+W^{[p]}\omega
\end{equation}
for any $\omega\in\Omega^p(M)$ where $\nabla$ (resp.\ $\nabla^*$) is the Levi-Civita connection \big(resp.\ its $L^2$-adjoint\big) on forms and \smash{$W^{[p]}$} is the curvature term. This last term is usually called the Weitzenb\"ock curvature operator and it defines a self-adjoint endomorphism acting on $p$-forms.
In the following, we will say that $W^{[p]}$ is bounded from below by $\gamma p(n-p)\in\mathbb{R}$ for a fixed $p=1,\dots,n$ if it satisfies
\begin{equation}\label{BCOperator}
\bigl\langle W^{[p]}\omega,\omega\bigr\rangle\geq \gamma p(n-p)|\omega|^2
\end{equation}
for all $\omega\in\Omega^p(M)$.
For $p=1$, this is equivalent to the fact that the Ricci curvature satisfies $\mathrm{Ric}\geq (n-1)\gamma g$.
From~\cite{GM}, this condition is satisfied for all $p$ if the curvature operator is bounded from below by $\gamma\in\mathbb{R}$.

Now integrating~(\ref{eq:weitzenboeckpforms}) on $M$ leads to the so-called {\it Reilly formula} on $p$-forms~\cite[Theorem~3]{RaulotSavo2011} which writes as
\begin{equation}\label{eq:Reillyformula}
\int_M\bigl(|{\rm d}\omega|^2+|\delta\omega|^2\bigr){\rm d}\mu=\int_M\bigl(|\nabla\omega|^2+\bigl\langle W^{[p]}\omega,\omega\bigr\rangle\bigr){\rm d}\mu+\int_{\partial M}\mathcal{B}(\omega,\omega){\rm d}\sigma
\end{equation}
for all $\omega\in\Omega^p(M)$ and where
\[
\mathcal{B}(\omega,\omega)=2\bigl\langle\nu\lrcorner\omega,\delta^{\partial M}(\iota^*\omega)\bigr\rangle+\bigl\langle S^{[p]}(\iota^*\omega),\iota^*\omega\bigr\rangle+(n-1)H|\nu\lrcorner\omega|^2-\bigl\langle S^{[p-1]}(\nu\lrcorner\omega),\nu\lrcorner\omega\bigr\rangle.
\]
As a direct consequence of this formula, we get the following estimate on $\lambda_{1,p}$.
\begin{Theorem}\label{GM_Dirichlet}
Let $(M^n,g)$ be a compact Riemannian manifold with boundary whose Weitzen\-b\"ock curvature operator $W^{[p]}$ is bounded from below by a positive constant $\gamma p(n-p)$ for some $1\leq p\leq [\frac{n}{2}]$.
Then we have $\lambda_{1,p}>\gamma p(n-p+1)$.
\end{Theorem}

\begin{proof} We first observe that if $\omega$ is a $p$-form on $M$ satisfying $\omega_{|_{\partial M}}=0$, the boundary term in the Reilly formula~(\ref{eq:Reillyformula}) vanishes and then, from~(\ref{BCOperator}), we get
\[
\int_M\bigl(|{\rm d}\omega|^2+|\delta\omega|^2\bigr){\rm d}\mu\geq\int_M\bigl(|\nabla\omega|^2+\gamma p(n-p)|\omega|^2\bigr){\rm d}\mu.
\]
On the other hand, with the help of the pointwise inequality
\[
|\nabla\alpha|^2\geq\frac{1}{p+1}|{\rm d}\alpha|^2+\frac{1}{n-p+1}|\delta\alpha|^2,
\]
which is true for any $p$-form $\alpha$ (see~\cite[Lemme 6.8]{GM}), we obtain
\[
\|{\rm d}\omega \|^2_{L^2(M)}+\|\delta \omega \|^2_{L^2(M)}\ge \frac{1}{n-p+1}\bigl(\|{\rm d}\omega \|^2_{L^2(M)}+\|\delta \omega \|^2_{L^2(M)}\bigr)+\gamma p(n-p)\|\omega \|^2_{L^2(M)}
\]
since \smash{$1\leq p\leq [\frac{n}{2}]$}. Note that equality occurs if $\omega$ is a conformal Killing form (see Lemma~\ref{CKFB} below for a precise definition and~\cite{Semmelmann02} for more general properties) and if $n=2p$ or if \smash{$1\leq p<[\frac{n}{2}]$} and ${\rm d}\omega=0$. We finally get
\begin{equation}\label{LaplaceBounds}
\|{\rm d}\omega\|^2_{L^2(M)}+\|\delta\omega \|^2_{L^2(M)} \geq \gamma p(n-p+1) \|\omega\|^2_{L^2(M)},
\end{equation}
which, by taking $\omega$ an eigenform associated to $\lambda_{1,p}$, leads to the desired estimate in a broad sense. Assume now that equality holds so that $\omega$ is a conformal Killing $p$-form which vanishes on the boundary. Then Lemma~\ref{CKFB} ensures that $\nabla_\nu\omega_{|\partial M}=0$ and this is impossible since $\omega$ is also a~$p$-eigenform for the Dirichlet Hodge Laplacian.
\end{proof}

In the previous proof, we considered conformal Killing $p$-forms $\omega$, that is $p$-forms satisfying the following equation
\begin{equation}\label{CKForms}
\nabla_X\omega=\frac{1}{p+1}X\lrcorner {\rm d}\omega-\frac{1}{n-p+1}X^\flat\wedge\delta\omega
\end{equation}
for all $X\in\Gamma(TM)$ and used the following result that we prove now.
\begin{Lemma}\label{CKFB}
On a $($not necessarily compact$)$ Riemannian manifold $(M^n,g)$ with boundary $\partial M$, a conformal Killing $p$-form with $1\leq p\leq n-1$ which vanishes on $\partial M$ satisfies $\nabla_\nu\omega_{|\partial M}=0$.
\end{Lemma}

\begin{proof} It is a direct consequence of~(\ref{CKForms}) that
\[
\iota^*(\nabla_\nu\omega)=\frac{1}{p+1}\nu\lrcorner {\rm d}\omega\qquad\text{and}\qquad\nu\lrcorner\nabla_\nu\omega=-\frac{1}{n-p+1}\iota^*(\delta\omega).
\]
On the other hand, since $\omega\in\Omega^p(M)$ vanishes on $\partial M$ it is straightforward from~(\ref{eq:iotastarnablanuomega}) and~(\ref{eq:nuintnablanuomega}) to compute that
\[
\iota^*(\nabla_\nu\omega)=\nu\lrcorner {\rm d}\omega\qquad\text{and}\qquad\nu\lrcorner\nabla_\nu\omega=-\iota^*(\delta\omega).
\]
Putting these relations together implies that
\[
\nu\lrcorner {\rm d}\omega=0 \qquad\text{and}\qquad\frac{n-p}{n-p+1}\iota^*(\delta\omega)=0,
\]
which, with the help of~(\ref{Zcharacterization}), allows us to conclude.
\end{proof}

Combining Theorem~\ref{GM_Dirichlet} with~(\ref{DirichletBounds}) leads to the following estimates for the buckling and clamped plate first eigenvalues.

\begin{Corollary}\label{BCDL}
Let $(M^n,g)$ be a compact Riemannian manifold with boundary whose Weitzen\-b\"ock curvature operator $W^{[p]}$ is bounded from below by $\gamma p(n-p)>0$ for some \smash{$1\leq p\leq [\frac{n}{2}]$}.
Then we~have
\[
\Lambda_{1,p}>\gamma p(n-p+1)\qquad\text{and}\qquad\Gamma_{1,p}>\gamma^2 p^2(n-p+1)^2.
\]
\end{Corollary}

\begin{Remarks}%\label{r:knownlowerboundsbucklingclampedplate}
\noindent\begin{enumerate}\itemsep=0pt
\item From these estimates, known lower bounds for the first eigenvalues of the buckling and clamped plate problems on functions~\cite{ChenChengWangXia2012} can be deduced.
Indeed, for $p=1$, the fact that~$W^{[1]}$ is bounded from below by $n-1$ is equivalent to the fact that the Ricci tensor of~$(M^n,g)$ satisfies $\mathrm{Ric}\geq (n-1)g$.
Then the first estimate in Corollary~\ref{BCDL} reads $\Lambda_{1,1}>n$ under this curvature assumption.
\item Independently, the third inequality of Theorem~\ref{BCP_Dirichlet}, which reads $\Lambda_{1,0}\geq \lambda_{1,1}$, combined with the inequality $\lambda_{1,1}>n$ from Theorem~\ref{GM_Dirichlet}, yields $\Lambda_{1,0}>n$.
This is exactly~\cite[Theorem~1.6]{ChenChengWangXia2012}.
Now putting together this estimate with the first inequality in Theorem~\ref{BCP_Dirichlet} gives $\Gamma_{1,0}>n\lambda_{1,0}$, which is precisely~\cite[Theorem 1.5]{ChenChengWangXia2012}.
\end{enumerate}
\end{Remarks}

\subsection{Domains in the unit sphere}%\label{ss:domainsinSn}

In this subsection, we derive an inequality which relates the first eigenvalues of the clamped plated problem for different degrees for domains of the $n$-dimensional unit sphere $\mathbb{S}^n$ carrying its metric of constant sectional curvature $1$.
This result is a consequence of a more general estimate for submanifolds isometrically immersed in the Euclidean space $\mathbb{R}^{n+m}$ which can be obtained using~\cite[Lemma 3.8]{EGHM}.
However, it is in general difficult to control all the terms which appear in these estimates so that we shall restrict ourselves to the case when $M$ is a domain of $\mathbb{S}^n$.
We~also obtain a similar result for the first buckling eigenvalue.
%when $n=2p$.
More precisely, we obtain the following.

\begin{Theorem}\label{t:domainSn}
Let $(M^n,g)$ be a domain in the unit sphere $\mathbb{S}^n$.
\begin{enumerate}\itemsep=0pt
\item[$1.$] For $1\leq p\leq [\frac{n}{2}]$, we have
\begin{equation}\label{eq:upperboundGammainSn}
p \Gamma_{1,p-1}+ (n-p) \Gamma_{1,p+1}< C_{n,p} \Gamma_{1,p},
\end{equation}
and
\begin{equation}\label{eq:upperboundLambdainSn}
\min( p\Lambda_{1,p-1}, (n-p)\Lambda_{1,p+1})< \frac{C_{n,p}}{2}\cdot\Lambda_{1,p},
\end{equation}
where \smash{$C_{n,p}=n+\frac{4+2(n-2p)^2}{p(n-p+1)}+\frac{n(n-2p)^2}{ p^2(n-p+1)^2}$}.
\item[$2.$] For $n=2p$, we have
\begin{equation}\label{eq:upperboundLambdainSnforn=2p}
\Lambda_{1,\frac{n}{2}-1} < \biggl(1+\frac{16}{n^2(n+2)}\biggr)\cdot \Lambda_{1,\frac{n}{2}}.
\end{equation}
\end{enumerate}
\end{Theorem}

\begin{proof} Let $\omega$ be a smooth $p$-form on $M$ which vanishes on $\partial M$ and whose normal covariant derivative is also zero on the boundary.
For any $1\leq i\leq n+1$, we consider the $(p-1)$-form \smash{$\partial_{x_i}^T\lrcorner\,\omega$} on $M$ where \smash{$\partial_{x_i}^T$} denotes the tangential part in $TM$ of the unit parallel vector field $\partial_{x_i}$ of~\smash{$\mathbb{R}^{n+1}$}. It is not difficult to check that
\[
\bigl(\partial_{x_i}^T\lrcorner\,\omega\bigr)_{|\partial M}=0\qquad\text{and}\qquad\nabla_{\nu}\bigl(\partial_{x_i}^T\lrcorner\omega\bigr)_{|\partial M}=0
\]
so that it can be used as a test-form in the variational characterizations of $\Gamma_{1,p-1}$ and $\Lambda_{1,p-1}$.

1. From~(\ref{eq:charcp}), we get
\[
\Gamma_{1,p-1}\int_M|\partial_{x_i}^T\lrcorner\,\omega|^2 {\rm d}\mu\leq\int_M|\Delta(\partial_{x_i}^T\lrcorner\,\omega)|^2{\rm d}\mu.
\]
Summing that inequality over $i$ and using the pointwise identity \smash{$\sum_{i=1}^{n+1}|\partial_{x_i}^T\lrcorner\,\omega|^2=p|\omega|^2$} (see, e.g.,~\cite[equation~(22)]{EGHM}), we obtain
\[
p\Gamma_{1,p-1}\int_M|\omega|^2 {\rm d}\mu\leq\int_M\sum_{i=1}^{n+1}|\Delta(\partial_{x_i}^T\lrcorner\,\omega)|^2{\rm d}\mu.
\]
Now the right-hand side of that inequality was computed for every $p$-form $\omega$ in the right-hand side of~\cite[equations~(31) and (33)]{EGHM} and we obtain
\begin{equation}\label{eq:inequalitysphere}
p\Gamma_{1,p-1} \|\omega\|^2_{L^2(M)} \leq 4\|\delta\omega\|^2_{L^2(M)}+p\|\Delta\omega+(2p-n)\omega\|^2_{L^2(M)}.
\end{equation}
Since inequality~\eqref{eq:inequalitysphere} is true for any $p$-eigenform $\omega$, we can apply it to the $(n-p)$-eigenform~$\star\omega$ to get
\[
(n-p) \Gamma_{1,n-p-1} \|\omega\|^2_{L^2(M)}\leq 4\|{\rm d}\omega\|^2_{L^2(M)}+(n-p)\|\Delta\omega+(n-2p)\omega\|^2_{L^2(M)}.
\]
Summing both inequalities and using the fact that $\Gamma_{1,n-p-1}=\Gamma_{1,p+1}$ yields
\begin{align*}%\label{eq:inesphere}
 \left(p \Gamma_{1,p-1}+ (n-p) \Gamma_{1,p+1}\right)\|\omega\|^2_{L^2(M)} \leq{}& n\|\Delta\omega\|^2_{L^2(M)}+n(2p-n)^2\|\omega\|^2_{L^2(M)}\nonumber\\
&{}{+}\,\bigl(4+2(n-2p)^2\bigr) (\Delta\omega,\omega)_{L^2(M)}\nonumber\\
\leq{}& n\|\Delta\omega\|^2_{L^2(M)}+n(2p-n)^2\|\omega\|^2_{L^2(M)}\nonumber\\
&{}{+}\,\bigl(4+2(n-2p)^2\bigr) \|\Delta\omega\|_{L^2(M)}\|\omega\|_{L^2(M)},\nonumber
\end{align*}
where we used the Cauchy--Schwarz inequality in the last inequality. On the other hand, since~$M$ is in the unit sphere, the estimate~(\ref{LaplaceBounds}) holds with $\gamma=1$, and therefore we get
\[
 (p \Gamma_{1,p-1}+ (n-p) \Gamma_{1,p+1} )\|\omega\|^2_{L^2(M)}\leq C_{n,p}\|\Delta\omega\|^2_{L^2(M)},
\]
which allows us to deduce~\eqref{eq:upperboundGammainSn} in the broad sense by taking $\omega$ to be an eigenform associated to $\Gamma_{1,p}$.
If~\eqref{eq:upperboundGammainSn} were an equality, then for any clamped-plate-eigenform $\omega$ associated to $\Gamma_{1,p}$, the $p$-form $\Delta\omega$ would be pointwise proportional to $\omega$ because of the equality in the Cauchy--Schwarz inequality. But, because of $\Delta^2\omega=\Gamma_{1,p}\omega$ on~$M$, we would have \smash{$\Delta\omega=\sqrt{\Gamma_{1,p}}\omega$} on~$M$, therefore $\omega$ would be a Dirichlet eigenform associated to the eigenvalue $\sqrt{\Gamma_{1,p}}$.
Again, because~$\omega$ satisfies $\nabla_\nu\omega_{|\partial M}=0$ along $\partial M$, the unique continuation property for elliptic second-order linear operators would imply that $\omega=0$ on $M$.
This would lead to a contradiction and shows that~\eqref{eq:upperboundGammainSn} must be strict.

Now, we prove~\eqref{eq:upperboundLambdainSn}. Using the variational characterization~(\ref{eq:charb}) of $\Lambda_{1,p-1}$ for $p=1,\dots,n$ gives
\[
\Lambda_{1,p-1}\int_{M}\bigl(\bigl|{\rm d}\bigl(\partial_{x_i}^T\lrcorner\omega\bigr)\bigr|^2+\bigl|\delta\bigl(\partial_{x_i}^T\lrcorner\omega\bigr)\bigr|^2\bigr){\rm d}\mu\leq \int_M\big|\Delta\bigl(\partial_{x_i}^T\lrcorner\omega\bigr)\big|^2{\rm d}\mu
\]
and summing as above on $i$ from $1$ to $n+1$ leads to
\[
\Lambda_{1,p-1}\sum_{i=1}^{n+1}\int_{M}\bigl(\bigl|{\rm d}\bigl(\partial_{x_i}^T\lrcorner\omega\bigr)\bigr|^2+\bigl|\delta\bigl(\partial_{x_i}^T\lrcorner\omega\bigr)\bigr|^2\bigr){\rm d}\mu\leq 4\|\delta\omega\|^2_{L^2(M)}+p\|\Delta\omega+(2p-n)\omega\|^2_{L^2(M)}.
\]
By the Cartan identity and~\cite[formula~(4.3)]{GueriniSavo2004} (see also~\cite[formula~(20)]{EGHM}), we have, for every $1\leq i\leq n+1$,
\begin{align*}
{\rm d}(\partial_{x_i}^T\lrcorner\omega)=\mathcal{L}_{\partial_{x_i}^T}\omega-\partial_{x_i}^T\lrcorner\,{\rm d}\omega
=\nabla_{\partial_{x_i}^T}\omega+\mathrm{I\!I}_{\partial^\perp_{x_i}}^{[p]}\omega-\partial_{x_i}^T\lrcorner\,{\rm d}\omega,
\end{align*}
where $\mathcal{L}_X\omega$ is the Lie derivative of $\omega$ in the $X$-direction and \smash{$\mathrm{I\!I}_{\partial^\perp_{x_i}}^{[p]}$} is the natural extension onto~$\Lambda^p T^*M$ of the pointwise endomorphism field \smash{$\mathrm{I\!I}_{\partial^\perp_{x_i}}$} of $TM$ defined for all $X,Y\in TM$ by
\[
\bigl\langle \mathrm{I\!I}_{\partial^\perp_{x_i}}(X),Y\bigr\rangle = \bigl\langle\mathrm{I\!I}(X,Y),\partial^\perp_{x_i}\bigr\rangle
\] and $\mathrm{I\!I}$ is the second fundamental form of $M$ in $\R^{n+1}$.
Here $\partial^\perp_{x_i}$ denotes the normal component of $\partial_{x_i}$ that is, $\partial^\perp_{x_i}=\partial_{x_i}-\partial^T_{x_i}$.
Note that, because $\mathrm{I\!I}_x=-g_x\otimes x$ for every $x\in\mathbb{S}^n$ and hence $x\in M$, we have \smash{$\mathrm{I\!I}_{\partial^\perp_{x_i}}=-x_i\cdot\mathrm{Id}$} at $x=(x_1,\dots,x_{n+1})$, so that \smash{$\mathrm{I\!I}_{\partial^\perp_{x_i}}^{[p]}=-px_i\cdot\mathrm{Id}$} at $x$.
It can be deduced that
\begin{align*}
\bigl|{\rm d} \bigl(\partial_{x_i}^T\lrcorner\omega\bigr)\bigr|^2 ={}&\Bigl|\nabla_{\partial_{x_i}^T}\omega+\mathrm{I\!I}_{\partial^\perp_{x_i}}^{[p]}\omega-\partial_{x_i}^T\lrcorner\,{\rm d}\omega\Bigr|^2\\
={}&\bigl|\nabla_{\partial_{x_i}^T}\omega\bigr|^2+\Bigl|\mathrm{I\!I}_{\partial^\perp_{x_i}}^{[p]}\omega\Bigr|^2+\bigl|\partial_{x_i}^T\lrcorner\,{\rm d}\omega\bigr|^2\\
&{}{+2}\Bigl\langle\nabla_{\partial_{x_i}^T}\omega,\mathrm{I\!I}_{\partial^\perp_{x_i}}^{[p]}\omega\Bigr\rangle-2\bigl\langle \nabla_{\partial_{x_i}^T}\omega,\partial_{x_i}^T\lrcorner\,{\rm d}\omega\bigr\rangle-2\Bigl\langle\mathrm{I\!I}_{\partial^\perp_{x_i}}^{[p]}\omega,\partial_{x_i}^T\lrcorner\,{\rm d}\omega\Bigr\rangle.
\end{align*}
Choosing any or\-tho\-nor\-mal basis $(e_j)_{1\leq j\leq n}$ of $T_x M$ for some $x\in M$, we have
\begin{align*}
\sum_{i=1}^n\bigl|\nabla_{\partial_{x_i}^T}\omega\bigr|^2&=\sum_{i=1}^{n+1}\sum_{j,k=1}^n\bigl\langle\partial^T_{ x_i},e_j\bigr\rangle\bigl\langle\partial^T_{ x_i},e_k\bigr\rangle\langle\nabla_{e_j}\omega,\nabla_{e_k}\omega\rangle\\
&=\sum_{i=1}^{n+1}\sum_{j,k=1}^n\langle\partial_{ x_i},e_j\rangle\langle\partial_{ x_i},e_k\rangle\langle\nabla_{e_j}\omega,\nabla_{e_k}\omega\rangle\\
&=\sum_{j,k=1}^n\Biggl(\sum_{i=1}^{n+1}\langle\partial_{ x_i},e_j\rangle\langle\partial_{ x_i},e_k\rangle\Biggr)\cdot\langle\nabla_{e_j}\omega,\nabla_{e_k}\omega\rangle\\
&=\sum_{j,k=1}^n\langle e_j,e_k\rangle\langle\nabla_{e_j}\omega,\nabla_{e_k}\omega\rangle
=\sum_{j=1}^n|\nabla_{e_j}\omega|^2
=|\nabla\omega|^2.
\end{align*}
As a second step, at every $x\in M$, because of $|x|=1$,
\[\sum_{i=1}^{n+1}\Bigl|\mathrm{I\!I}_{\partial^\perp_{x_i}}^{[p]}\omega\Bigr|^2=p^2\sum_{i=1}^{n+1}x_i^2|\omega|^2=p^2|x|^2|\omega|^2=p^2|\omega|^2.\]
By~\cite[equation~(22)]{EGHM},
\[\sum_{i=1}^{n+1}\bigl|\partial_{x_i}^T\lrcorner\,{\rm d}\omega\bigr|^2=(p+1)|{\rm d}\omega|^2.\]
Moreover,
\begin{align*}
\sum_{i=1}^{n+1}2\Bigl\langle\nabla_{\partial_{x_i}^T}\omega,\mathrm{I\!I}_{\partial^\perp_{x_i}}^{[p]}\omega\Bigr\rangle& =p\sum_{i=1}^{n+1}2x_i\bigl\langle\nabla_{\partial_{x_i}^T}\omega,\omega\bigr\rangle
=p\sum_{i=1}^{n+1}x_i\partial_{x_i}^T\bigl(|\omega|^2\bigr)
=px^T\bigl(|\omega|^2\bigr)
=0
\end{align*}
because of $x^T=0$ for every $x\in\mathbb{S}^n$.
For the same reason,
\begin{align*}
\sum_{i=1}^{n+1}\Bigl\langle\mathrm{I\!I}_{\partial^\perp_{x_i}}^{[p]}\omega,\partial_{x_i}^T\lrcorner\,{\rm d}\omega\Bigr\rangle&=p\sum_{i=1}^{n+1}x_i\bigl\langle\omega,\partial_{x_i}^T\lrcorner\,{\rm d}\omega\bigr\rangle
=p\bigl\langle\omega,x^T\lrcorner\,{\rm d}\omega\bigr\rangle
=0.
\end{align*}
Applying the same computational method as above, we also have, in any pointwise or\-tho\-nor\-mal basis $(e_j)_{1\leq j\leq n}$ of $TM$,
\begin{align*}
\sum_{i=1}^{n+1}\langle \nabla_{\partial_{x_i}^T}\omega,\partial_{x_i}^T\lrcorner\,{\rm d}\omega\rangle&=\sum_{j=1}^n\langle \nabla_{e_j}\omega,e_j\lrcorner\,{\rm d}\omega\rangle
=\sum_{j=1}^n\bigl\langle e_j^\flat\wedge\nabla_{e_j}\omega,{\rm d}\omega\bigr\rangle
=|{\rm d}\omega|^2.
\end{align*}
On the whole, we obtain
\begin{align*}
\sum_{i=1}^{n+1}\bigl|{\rm d}(\partial_{x_i}^T\lrcorner\,\omega)\bigr|^2&=|\nabla\omega|^2+p^2|\omega|^2+(p+1)|{\rm d}\omega|^2-2|{\rm d}\omega|^2\\
&=|\nabla\omega|^2+p^2|\omega|^2+(p-1)|{\rm d}\omega|^2.
\end{align*}
Therefore,
\[
\sum_{i=1}^{n+1}\bigl\|{\rm d}\bigl(\partial_{x_i}^T\lrcorner\,\omega\bigr)\bigr\|_{L^2(M)}^2= \|\nabla\omega\|_{L^2(M)}^2+p^2\|\omega\|_{L^2(M)}^2+(p-1)\|{\rm d} \omega\|_{L^2(M)}^2.
\]
Using the Weitzenb\"ock formula~(\ref{eq:weitzenboeckpforms}) with $W^{[p]}=p(n-p)\cdot\mathrm{Id}$ and \smash{$\omega_{|_{\partial M}}=\nabla_\nu\omega_{|_{\partial M}}=0$}
leads to
\begin{align*}
\sum_{i=1}^{n+1}\bigl\|{\rm d}\bigl(\partial_{x_i}^T\lrcorner\,\omega\bigr)\bigr\|_{L^2(M)}^2={}&(\nabla^*\nabla\omega,\omega)_{L^2(M)} +p^2\|\omega\|_{L^2(M)}^2+(p-1)\|{\rm d}\omega\|_{L^2(M)}^2\\
={}&(\Delta\omega,\omega)_{L^2(M)}-\bigl(W^{[p]}\omega,\omega\bigr)_{L^2(M)}+p^2\|\omega\|_{L^2(M)}^2\\
&{}{+}\,(p-1)\|{\rm d}\omega\|_{L^2(M)}^2\\
={}&(\Delta\omega,\omega)_{L^2(M)}+\bigl(p^2-p(n-p)\bigr)\|\omega\|_{L^2(M)}^2+(p-1)\|{\rm d}\omega\|_{L^2(M)}^2\\
={}&p(2p-n)\|\omega\|_{L^2(M)}^2+(p-1)\|{\rm d}\omega\|_{L^2(M)}^2+(\Delta\omega,\omega)_{L^2(M)}.
\end{align*}
On the other hand, for every $i\in\{1,\dots,n+1\}$, we can compute, using a local or\-tho\-nor\-mal basis $(e_j)_{1\leq j\leq n}$ of $TM$ as well as $\nabla_X\partial_{x_i}^T=-x_i X$ for every $X\in T_xM$,
\begin{align*}
\delta\bigl(\partial_{x_i}^T\lrcorner\,\omega\bigr)&=-\sum_{j=1}^ne_j\lrcorner\,\nabla_{e_j}\bigl(\partial_{x_i}^T\lrcorner\,\omega\bigr)\\
&=-\sum_{j=1}^ne_j\lrcorner\bigl(\bigl(\nabla_{e_j}\partial_{x_i}^T\bigr)\lrcorner\,\omega+\partial_{x_i}^T\lrcorner\,\nabla_{e_j}\omega\bigr)\\
&=-\sum_{j,k=1}^n\bigl\langle\nabla_{e_j}\partial_{x_i}^T,e_k\bigr\rangle e_j\lrcorner\,e_k\lrcorner\,\omega-\partial_{x_i}^T\lrcorner\,\delta\omega\\
&=\underbrace{\sum_{j,k=1}^nx_i\langle e_j,e_k\rangle e_j\lrcorner\,e_k\lrcorner\,\omega}_{0}-\partial_{x_i}^T\lrcorner\,\delta\omega\\
&=-\partial_{x_i}^T\lrcorner\,\delta\omega,
\end{align*}
where we used the skew-symmetry of $\langle e_j,e_k\rangle e_j\lrcorner\,e_k\lrcorner\,\omega$ in $(j,k)$.
As a consequence, again by~\cite[equation~(22)]{EGHM},
\begin{align*}
\sum_{i=1}^{n+1}\bigl|\delta(\partial_{x_i}^T\lrcorner\omega)\bigr|^2=\sum_{i=1}^{n+1}\bigl|\partial_{x_i}^T\lrcorner\,\delta\omega\bigr|^2=(p-1)|\delta\omega|^2,
\end{align*}
from which \smash{$ {\sum_{i=1}^{n+1}\bigl\|\delta\bigl(\partial_{x_i}^T\lrcorner\omega\bigr)\bigr\|_{L^2(M)}^2=(p-1)\|\delta\omega\|_{L^2(M)}^2}$} follows.
Finally, we deduce that
\begin{gather*}
 \Lambda_{1,p-1}\bigl(p(2p-n)\|\omega\|_{L^2(M)}^2+(p-1)\bigl(\|{\rm d}\omega\|_{L^2(M)}^2+\|\delta\omega\|_{L^2(M)}^2\bigr)+(\Delta\omega,\omega)_{L^2(M)}\bigr)\\
\qquad\leq4\|\delta\omega\|^2_{L^2(M)}+p\|\Delta\omega+(2p-n)\omega\|^2_{L^2(M)},
\end{gather*}
that is, using \smash{$\|{\rm d}\omega\|_{L^2(M)}^2+\|\delta\omega\|_{L^2(M)}^2=(\Delta\omega,\omega)_{L^2(M)}$},
\begin{gather*}
%\label{buckling_shpere_pminus1}
p\Lambda_{1,p-1}\bigl((\Delta\omega,\omega)_{L^2(M)}+(2p-n)\|\omega\|^2_{L^2(M)}\bigr)\leq 4\|\delta\omega\|^2_{L^2(M)}+p\|\Delta\omega+(2p-n)\omega\|^2_{L^2(M)}.
\end{gather*}
Notice that the left-hand side of that last identity must be positive since it is actually
\[
\sum_{i=1}^{n+1}\bigl\|{\rm d}\bigl(\partial_{x_i}^T\lrcorner\omega\bigr)\bigr\|_{L^2(M)}^2+\bigl\|\delta\bigl(\partial_{x_i}^T\lrcorner\omega\bigr)\bigr\|_{L^2(M)}^2.
\]
Replacing $\omega$ by $\star\omega$ and using the Hodge symmetry of the buckling, eigenvalues ensure that
\begin{gather*}
\begin{split}
&(n-p)\Lambda_{1,p+1}\bigl((\Delta\omega,\omega)_{L^2(M)}+(n-2p)\|\omega\|^2_{L^2(M)}\bigr)\\
& \qquad\leq 4\|{\rm d}\omega\|^2_{L^2(M)}+(n-p)\|\Delta\omega+(n-2p)\omega\|^2_{L^2(M)}.
\end{split}
\end{gather*}
Adding those two inequalities, we obtain
\begin{gather}
\nonumber 2\min(p\Lambda_{1,p-1}, (n-p)\Lambda_{1,p+1})\cdot (\Delta\omega,\omega)_{L^2(M)}
\\
\nonumber\qquad{}\le{}4\bigl(\|{\rm d}\omega\|^2_{L^2(M)}+\|\delta\omega\|^2_{L^2(M)}\bigr)+p\|\Delta\omega+(2p-n)\omega\|^2_{L^2(M)}\\
\nonumber\qquad\quad{}+(n-p)\|\Delta\omega+(n-2p)\omega\|^2_{L^2(M)}
\\
\nonumber\qquad{}=4 (\Delta\omega,\omega )_{L^2(M)}+n\|\Delta\omega\|_{L^2(M)}^2 +2 (p(2p-n)+(n-p)(n-2p) ) (\Delta\omega,\omega )_{L^2(M)}\\
\nonumber\qquad\quad{}+\bigl(p(2p-n)^2+(n-p)(n-2p)^2\bigr)\|\omega\|_{L^2(M)}^2\\
\qquad{}=n\|\Delta\omega\|^2_{L^2(M)}+n(n-2p)^2\|\omega\|^2_{L^2(M)}+\bigl(4+2(n-2p)^2\bigr)(\Delta\omega,\omega)_{L^2(M)}.\label{buckling_sphere}
\end{gather}
But~\eqref{LaplaceBounds} with $\gamma=1$ gives
\[\|\omega\|_{L^2(M)} \le \frac{1}{p(n-p+1)} \|\Delta \omega\|_{L^2(M)}.\]
Substituting that estimate in~\eqref{buckling_sphere} yields the following inequality:
\[2 \min(p\Lambda_{1,p-1}, (n-p)\Lambda_{1,p+1})\cdot\bigl(\|{\rm d}\omega\|_{L^2(M)}^2+\|\delta\omega\|_{L^2(M)}^2\bigr) \le C_{n,p} \|\Delta \omega\|^2_{L^2(M)},\]
where $C_{n,p}$ is the constant defined in the first statement of Theorem~\ref{t:domainSn}.
We deduce inequality~(\ref{eq:upperboundLambdainSn}) by taking $\omega$ as an eigenform of the buckling eigenvalue problem associated to $\Lambda_{1,p}$. Similarly to the first case, the equality cannot occur in~(\ref{eq:upperboundLambdainSn}).

2. If $n=2p$, then~(\ref{eq:upperboundLambdainSn}) becomes, because of $\Lambda_{1,p-1}=\Lambda_{1,n-p+1}=\Lambda_{1,p+1}$,
\[\Lambda_{1,\frac{n}{2}-1}<\frac{C_{n,\frac{n}{2}}}{n}\cdot\Lambda_{1,\frac{n}{2}},\]
which is~(\ref{eq:upperboundLambdainSnforn=2p}).
\end{proof}

Notice that, by the assumption $p\leq n-p$, inequality~(\ref{eq:upperboundLambdainSn}) implies that
\[\min(\Lambda_{1,p-1},\Lambda_{1,p+1})< \frac{C_{n,p}}{2p}\cdot\Lambda_{1,p}.\]

\subsection*{Acknowledgements}
This collaboration is supported by the International Emerging Action ``HOPF'' (Higher Order boundary value Problems for differential Forms) of the French CNRS, which the authors would like to thank.
We also acknowledge the support of Universit\'e de Rouen Normandie, Universit\'e de Lorraine, and Universit\'e Libanaise.
{\sl Last but not least}, we thank the referees for their careful reading and constructive comments as well as the editorial team for its efficient work.

\pdfbookmark[1]{References}{ref}
\LastPageEnding

\end{document}